\documentclass{article}

\usepackage {amsmath}
\usepackage {hyperref}
\usepackage {amsfonts}
\usepackage {amsthm}
\usepackage[flushmargin]{footmisc}
\usepackage[parfill]{parskip}
\usepackage {amssymb}
\usepackage {fullpage}
\usepackage {mathrsfs}
\usepackage {amscd}
\usepackage[all,cmtip]{xy}
\DeclareMathAlphabet{\mathpzc}{OT1}{pzc}{m}{it}

\title{Shintani descent for algebraic groups and almost characters of unipotent groups}
\author{Tanmay Deshpande}
\date{}

\newtheorem {thm} {Theorem} [section]
\newtheorem {prop} [thm] {Proposition}
\newtheorem {conj} [thm] {Conjecture}
\newtheorem {lem} [thm] {Lemma}
\newtheorem {cor} [thm] {Corollary}

\theoremstyle{definition}
\newtheorem {defn} [thm] {Definition}

\newtheorem {prob} [thm]  {Problem}

\newtheorem {rk} [thm]  {Remark}
\newtheorem {ex} [thm] {Example}

\newcommand{\beq}{\begin{equation}}
\newcommand{\eeq}{\end{equation}}
\newcommand{\bthm}{\begin {thm}}
\newcommand{\ethm}{\end {thm}}
\newcommand{\bprop}{\begin {prop}}
\newcommand{\eprop}{\end {prop}}
\newcommand{\bprob}{\begin {prob}}
\newcommand{\eprob}{\end {prob}}
\newcommand{\bcor}{\begin {cor}}
\newcommand{\ecor}{\end {cor}}
\newcommand{\blem}{\begin{lem}}
\newcommand{\elem}{\end{lem}}
\newcommand{\bdefn}{\begin{defn}}
\newcommand{\edefn}{\end{defn}}
\newcommand{\brk}{\begin{rk}}
\newcommand{\erk}{\end{rk}}

\newcommand{\xto}{\xrightarrow}

\renewcommand {\bar} {\overline}
\newcommand{\bpf}{\begin{proof}}
\newcommand{\epf}{\end{proof}}
\newcommand{\bex}{\begin{ex}}
\newcommand{\eex}{\end{ex}}
\newcommand{\rar}[1]{\stackrel{#1}{\longrightarrow}}
\newcommand{\f}{\mathbb}

\newcommand{\<}{\langle}
\renewcommand{\>}{\rangle}
\newcommand{\cycF}{\< F\rangle}

\newcommand{\h}{\operatorname}

\newcommand{\End}{\operatorname{End}}
\newcommand{\Sh}{\operatorname{Sh}}
\newcommand{\ad}{\operatorname{ad}}

\renewcommand{\L}{\mathcal{L}}

\newcommand{\Fq} {\mathbb{F}_q}

\newcommand{\Fqcl} {\overline{\mathbb{F}}_q}
\newcommand{\Q} {\mathbb{Q}}

\newcommand{\F} {\mathbb{F}}

\renewcommand{\t}{\widetilde}
\renewcommand{\phi} {\varphi}

\renewcommand{\k} {\mathtt{k}}

\newcommand{\M} {\mathcal{M}}

\newcommand{\D} {\mathscr{D}}

\newcommand{\uuPic} {\underline{\underline{\h{Pic}}}}
\newcommand{\uPic} {\underline{\h{Pic}}}
\newcommand{\Pic} {{\h{Pic}}}
\newcommand{\uEqBr} {\underline{\h{EqBr}}}
\newcommand{\EqBr} {{\h{EqBr}}}
\renewcommand{\Vec}{\h{Vec}}
\newcommand{\Aut}{\h{Aut}}

\newcommand{\Stab}{\h{Stab}}

\newcommand{\Qlcl} {\overline{\mathbb{Q}}_{\ell}}

\newcommand{\Rab} {R_{\g_1,\g_2}}
\newcommand{\Rba} {R_{\g_2,\g_1}}
\newcommand{\Rif} {R_{\id_G,F}}

\newcommand{\noin}{\noindent}

\newcommand{\g} {{\gamma}}

\newcommand{\id}{\operatorname{id}}

\newcommand{\Meg}{\M_{G,e}}

\newcommand{\bit}{\begin{itemize}}
\newcommand{\eit}{\end{itemize}}

\newcommand{\Rep}{\h{Rep}}
\newcommand{\Irrep}{\h{Irrep}}

\newcommand{\DG}{\D_G(G)}

\newcommand{\DgaG}{\D_G^{\g_1}(G)}

\newcommand{\DGab}{\D_G^{\g_1}(G)^{\g_2}}

\newcommand{\FunG}{\Fun([G],F)}

\newcommand{\E}{\mathcal{E}}

\newcommand{\Z}{\mathbb{Z}}

\newcommand{\tr}{\h{tr}}
\newcommand{\Fun}{\h{Fun}}

\newcommand{\bconj}{\begin{conj}}
\newcommand{\econj}{\end{conj}}
\begin{document}
\maketitle

\begin{abstract}
\noin In this paper, we extend the notion of Shintani descent to general (possibly disconnected) algebraic groups defined over a finite field $\Fq$. For this, it is essential to treat all the pure inner $\Fq$-rational forms of the algebraic group at the same time. We prove that the notion of almost characters (introduced by T. Shoji using Shintani descent) is well defined for any neutrally unipotent algebraic group, i.e. an algebraic group whose neutral connected component is a unipotent group. We also prove that these almost characters coincide with the ``trace  of Frobenius'' functions associated with Frobenius-stable character sheaves on neutrally unipotent groups. In the course of the proof, we also prove that the modular categories that arise from Boyarchenko-Drinfeld's theory of character sheaves on neutrally unipotent are in fact positive integral confirming a conjecture due to Drinfeld.
\end{abstract}

\tableofcontents

\section*{Acknowledgments}
I am grateful to Toshiaki Shoji for some illuminating conversations which led to this work. This work was supported by World Premier Institute Research Center Initiative (WPI), MEXT, Japan.

\section{Introduction}
Let us fix a prime number $p$, $q$ a power of $p$ and let us set $\k=\Fqcl$. Let $G$ be an algebraic group over $\k$ equipped with an $\Fq$-structure defined by a Frobenius $F:G\to G$. Shintani descent on $G$ compares the irreducible representations of the finite groups $G(\F_{q^m})$ and $G(\Fq)$ where $m$ is a positive integer. We begin by recalling the notion of Shintani descent for connected algebraic groups from \cite{S}. 

\subsection{Shintani descent for connected algebraic groups}
We refer to \cite{S} for the details of the constructions in this introduction. On the other hand, all these constructions will be carried out in detail in a more general setting in the subsequent sections of this paper.

Suppose that $G$ is connected. Shintani descent for $G$ was studied by Shoji in \cite{S} extending the work \cite{Shi} of Shintani  for general linear groups. Let $m$ be any positive integer and consider the $m$-th power of the Frobenius $$F^m:G\to G.$$ This is a Frobenius map for an $\F_{q^m}$-rational structure on $G$. The Frobenius $F$ restricts to a group automorphism $F:G^{F^m}\to G^{F^m}$ of order $m$ and hence it acts on the finite set $\Irrep (G^{F^m})$ of irreducible characters of $G^{F^m}$. Shintani descent gives us a way to compare the two sets $\Irrep(G^{F^m})^F$ and $\Irrep(G^F)$.

More precisely, consider the action of $G^{F^m}$ on itself by $F$-twisted conjugation defined by 
\beq\label{twistconj} g:h\mapsto ghF(g^{-1}).\eeq
Let us denote the set of orbits of this action by $G^{F^m}/{\sim_F}$. Let $\Fun(G^{F^m}/{\sim_F})$ denote the space of $F$-conjugation invariant $\Qlcl$-valued functions on $G^{F^m}$ where $\ell\neq p$ is a prime number. 

\brk
In this paper we will always considers functions and characters with values in $\Qlcl$ since we will use the theory of $\Qlcl$-sheaves and complexes. Once we fix an isomorphism $\Qlcl\cong \f{C}$, we can talk of ``complex conjugation'' on the field $\Qlcl$ and we can define the usual Hermitian inner product on the function spaces \linebreak $\Fun(G^{F^m}/\sim_F)\subset\Fun(G^{F^m})$. The choice of the isomorphism  $\Qlcl\cong \f{C}$ does not matter because in this paper we will only need to consider ``complex conjugation'' on the subfield $\Q^{ab}\subset \Qlcl$ where ``complex conjugation'' can be canonically defined.
\erk

Now given $\chi\in \Irrep(G^{F^m})^F$, we can define\footnote{Extend $\chi$ to an irreducible character $\chi'$ of the semidirect product $G^{F^m}\langle F\rangle$, where $\cycF$ is the cyclic group of order $m$ generated by $F$. Then set $\t\chi(g)=\chi'(gF)$ for $g\in G^{F^m}$. Note that there is a choice involved in extending $\chi$ to $\chi'$ and the different choices differ by scaling by $m$-th roots of unity on the coset $G^{F^m}\cdot F\subset G^{F^m}\<F\>$.} a function $\t\chi\in \Fun(G^{F^m}/{\sim_F})$ (well defined up to scaling by $m$-th roots of unity). The set $\{\t\chi|\chi\in \Irrep(G^{F^m})^F\}$ forms an orthonormal basis of the Hermitian inner product space $\Fun(G^{F^m}/{\sim_F})$.

Using Lang's theorem we can define the norm map which is a bijection between the set $G^{F^m}/{\sim_F}$ and the set $G^F/\sim$ of conjugacy classes in $G^F$:
\beq
N_m:\left(G^{F^m}/{\sim_F}\right)\rar{\cong}\left(G^{F}/{\sim}\right).
\eeq

Hence we get an isomorphism of function spaces
\beq\label{normmap}
{N_m^{-1}}^*:\Fun(G^{F^m}/{\sim_F})\xto{\cong}\Fun(G^{F}/{\sim}).
\eeq

\bdefn 
For a positive integer $m$, Shintani descent is defined to be the map
\beq\label{shincon}
\Sh_m:\Irrep(G^{F^m})^F \hookrightarrow \Fun(G^F/\sim)
\eeq
defined by 
\beq
\Irrep(G^{F^m})^F\ni \chi \mapsto {N_m^{-1}}^*(\t\chi)\in \Fun(G^F/\sim).
\eeq
\edefn
\brk
The map $\Sh_m$ is only well defined up to scaling by $m$-th roots of unity since this is so for the assignment $\Irrep(G^{F^m})^F\ni\chi\mapsto \t\chi\in \Fun(G^{F^m}/\sim_F)$.
\erk

The image of $\Sh_m$ forms an orthonormal basis of the Hermitian space $\Fun(G^F/\sim)$ which we call the $m$-th Shintani basis. The set $\Irrep(G^F)$ of irreducible characters of $G^F$ gives another orthonormal basis. 

\brk\label{r:shintanigln}
In general these two bases of $\Fun(G^F/\sim)$ are not equal (even for $m=1$, see \S\ref{s:twistoper}). In fact, describing the relationship between these two bases is one of the objectives of this paper. However, if $G=GL_n$ as in the original case studied by Shintani, then these bases are in fact equal and Shintani descent defines an explicit bijection between the sets $\Irrep(GL_n(\F_{q^m}))^F$ and $\Irrep(GL_n(\Fq))$. 
\erk

In \cite{S} Shoji proved that if $G$ is connected and satisfies certain further conditions, then the $m$-th Shintani basis of $\Fun(G^F/\sim)$ is (up to scaling by roots of unity) independent  of $m$ if $m$ is \emph{sufficiently divisible}. The elements of this common Shintani basis (when it exists) for sufficiently divisible $m$ are known as the almost characters of $G^F$. Shoji proved that the almost characters are eigenvectors for the twisting operator defined by setting $m=1$ in (\ref{normmap}): 
\beq\label{twistcon}
N_1^*:\Fun(G^F/\sim)\to \Fun(G^F/\sim).
\eeq
 Moreover, if $G$ is reductive and satisfies some further conditions (for example say if $G$ has connected center or $G$ is special linear and if $p$ is not too small) then Shoji has proved that the almost characters of $G^F$ coincide with the ``trace of Frobenius'' functions associated with $F$-stable character sheaves on $G$. Our main result, Theorem \ref{mainthm} below, says that analogues of all these results are true for neutrally unipotent groups $G$, namely groups such that their neutral connected components $G^\circ$ are unipotent algebraic groups.

\subsection{Shintani descent in general and pure inner forms}
In \cite{S}, an analogous descent construction is also carried out for disconnected groups under some additional conditions. In \S\ref{s:sdgag}, we define Shintani descent in general. The main difficulty in treating the disconnected case is the fact that the Lang isogeny for disconnected groups may not be surjective. However, this can be overcome if we consider not just the Frobenius $F$ that we start out with, but also all pure inner forms of $F$. We will now briefly describe the main results of this paper.

Let $G$ be any (possibly disconnected) algebraic group over $\k$. Let $F:G\to G$ be a Frobenius map defining an $\Fq$-rational structure on $G$. Consider the $F$-twisted conjugation action of $G$ on itself as described by (\ref{twistconj}). Let $H^1(F,G)$ denote the set of $F$-twisted conjugacy classes in $G$. By Lang's Theorem we have a natural bijection $H^1(F,G)=H^1(F,\Pi_0(G))$. In particular $H^1(F,G)$ is a finite set and each orbit of the $F$-twisted conjugation action is a union of  some of the connected components of the group $G$. 

For each $g\in G$, we can define a new Frobenius map $gF:=\ad(g)\circ F:G\to G$ which we call an inner form of $F$. The isomorphism class of  the $\Fq$-rational form of $G$ defined by the Frobenius $gF$ only depends on the $F$-twisted conjugacy class of $g$ since we have a commutative diagram
$$\xymatrixcolsep{7pc}\xymatrix{
G\ar[r]^-{gF}\ar[d]_{ad(h)} & G\ar[d]^{ad(h)}\\
G\ar[r]^-{{hgF(h^{-1})}F} & G\\
}$$
for each $h,g\in G$. We are interested in the representation theory of $G^F$. However, it is more natural to consider all pure inner forms $G^{\ad(g)\circ F}$ as $g\in G$ goes over a set of representatives of $H^1(F,G)$.  

\bdefn\label{d:irrepfun}
We define the set 
\beq
\Irrep(G,F):=\coprod\limits_{\< g\>\in H^1(F,G)} \Irrep(G^{\ad(g)\circ F})=\coprod\limits_{\< g\>\in H^1(F,G)} \Irrep(G^{gF})
\eeq
and the commutative $\Qlcl$-algebra (under convolution of functions)
\beq
\Fun([G],F):=\prod\limits_{\< g\>\in H^1(F,G)} \Fun(G^{gF}/\sim).
\eeq
We equip $\FunG$ with the standard Hermitian inner product. By taking the characters of irreducible representations, we can consider the set $\Irrep(G,F)$ as an orthonormal  basis of $\Fun([G],F)$.
\edefn

In \S\ref{s:defsd}, for each positive integer $m$ we will define a Shintani descent map (well defined up to scaling by $m$-th roots of unity)
\beq
\Sh_{m}:\Irrep (G,F^m)^F\hookrightarrow \Fun([G], F)
\eeq
generalizing (\ref{shincon}).  The image of $\Sh_m$ inside $\Fun([G],F)$ is an orthonormal basis which we again call the $m$-th Shintani basis (which is well defined up to scaling by $m$-th roots of unity).

Also, just as in the connected case (see (\ref{twistcon})), we will construct a (unitary) twisting operator
\beq\label{twistintro}
\Theta=N_1^*:\FunG\to \FunG.
\eeq
The map $N_1$ that we will define in \S\ref{s:defsd} permutes the set of rational conjugacy classes of all the pure inner forms within each geometric conjugacy class.

\brk
Each irreducible character $\chi\in \Irrep(G,F)\subset \Fun([G],F)$ is supported on exactly one pure inner form $G^{gF}$. However, if $\chi\in \Irrep(G,F^m)^F$, then $\Sh_m(\chi)\in \FunG$ can be nonzero on multiple pure inner forms. Also the map $N_1$ can map a rational conjugacy class of one inner form to a rational conjugacy class of a different inner form. From these considerations we see that it is essential to consider all pure inner forms $G^{gF}$ at the same time.
\erk

We now introduce the category $\DG$. This is the triangulated braided monoidal category of conjugation equivariant $\Qlcl$-complexes on $G$. The monoidal structure comes from the convolution with compact support of $\Qlcl$-complexes. For each object $C\in \DG$, we have functorial automorphisms $\theta_C:C\to C$ known as twists. We refer to \cite{BD} for more about the structure of the category $\DG$. We recall from \cite[\S2.4.8]{De} that given an object $C\in \DG$ and an isomorphism $\psi:F^*C\to C$, we have the associated ``trace of Frobenius'' function $T_{C,\psi}\in \FunG$.

We can now state our main result. Let $G$ be an algebraic group such that its neutral connected component $G^\circ$ is a unipotent group. In this case, the notion of character sheaves on $G$ has been defined\footnote{To be precise, in {\cite{BD}} character sheaves are only studied for unipotent groups. However everything can be readily extended to the neutrally unipotent case.} in \cite{BD}. We let $CS(G)$ denote the set of character sheaves on $G$. These are (the isomorphism classes of) certain special objects in the category $\DG$.
\bthm\label{mainthm}
Let $G$ be a neutrally unipotent group over $\k$ as above. Let $F:G\to G$ be an $\Fq$-Frobenius. Then we have the following:\\
(i) There exists a positive integer $m_0$ such that if $m$ is divisible by $m_0$, then the $m$-th Shintani basis of $\FunG$ is (up to scaling by roots of unity) independent of $m$. We define almost characters to be the elements of this common Shintani basis.\\
(ii) The almost characters as defined above are eigenvectors for the twisting operator $\Theta$ (see (\ref{twistintro})).\\
(iii) Recall that the Frobenius $F$ induces a permutation of the set $CS(G)$. The mapping which takes an $F$-stable character sheaf to its associated trace of Frobenius function in $\FunG$ defines a bijection from the set $CS(G)^F$ to the set of almost characters.
\ethm

\brk
We conjecture that in fact the above  results should hold for all algebraic groups. In particular we expect that there should be an interesting theory of character sheaves on all algebraic groups.  
\erk

In \S\ref{s:prelcon} we describe some preliminary constructions. In particular, for each $m\in \Z$ we introduce the $\Qlcl$-linear triangulated category $\D^{F^m}_G(G)$ which encodes the representation theory of $G^{F^m}$ and all its pure inner forms and define associative convolution products $\D^{F^{m_1}}_G(G)\times \D^{F^{m_2}}_G(G)\rar{}\D^{F^{m_1+m_2}}_G(G)$ which satisfy some crossed braiding relations. We construct analogues of the sheaf-function correspondence in certain general settings.

Then we use the various constructions from \S\ref{s:prelcon} to define Shintani descent for general algebraic groups in \S\ref{s:defsd}. In \S\ref{s:twistoper} we study the special case of Shintani descent for $m=1$, namely the twisting operator introduced in (\ref{twistintro}). In \S\ref{s:twistandsf} we prove that the trace functions associated with $F$-stable simple objects in $\DG$ are always eigenvectors for the twisting operator $\Theta$ from (\ref{twistintro}).

In \S\ref{s:twisttrace} we define the notion of twists in the categories $\D^{F^m}_G(G)$ and study their properties. We prove that these twists are compatible with the traces in $\D^{F^{m}}_G(G)$ that were defined in \cite[\S2.4.5]{De}.

In \S\ref{s:ipf} we express the inner products between elements of the $m$-th Shintani basis and irreducible characters in terms of the traces and convolution structure of the categories $\D^{F^m}_G(G)$. Note that these inner products are the entries of the unitary matrix which relates the $m$-th Shintani basis of $\FunG$ to the orthonormal basis of $\FunG$ formed by the irreducible characters of all the pure inner forms.

From \S\ref{s:pfmain} onwards we restrict to the case of neutrally unipotent groups. In \S\ref{s:pfmain} we recall some relevant facts from the theory of character sheaves on unipotent groups. In particular we recall the notion of $\f{L}$-packets and prove that Shintani descent respects the $\f{L}$-packet decompositions. We state a refined version of Theorem \ref{mainthm}. 

In \S\ref{s:catprel} we study some properties of the modular categories and their module categories that arise from the theory of character sheaves on neutrally unipotent groups. In particular we will prove (see Theorem \ref{t:posint}) that the modular categories that arise from this theory are positive integral.

Finally in \S\ref{s:pfref} we complete the proof of the refined version of our main theorem.

\section{Preliminary constructions}\label{s:prelcon}
In this section we will describe some preliminary constructions which will help us to define and study Shintani descent for general algebraic groups over $\k$ equipped with a Frobenius map. In this paper, by passing to the perfectizations, we will always assume that all our groups are in fact perfect quasi-algebraic groups over $\k$ even though we may not mention this explicitly. We refer to \cite[\S1.9]{BD} for more about this convention. With this convention, the Frobenius maps in fact become automorphisms.

\subsection{Twisted conjugation actions}\label{s:tconjact}
Let $G$ be a (perfect quasi-) algebraic group over $\k$ and let $\g:G\to G$ be an automorphism. We define the $\g$-conjugation action of $G$ on itself by $g:h\mapsto gh\g(g)^{-1}$. Equivalently, we can think of this action as the action of $G$ by conjugation on the coset $G\g$ in the semidirect product $\t G:=G\<\g\>=G\rtimes \Z$ under the natural identification of $G$ with $G\g$. We let $\D_G^\g(G)$ denote the equivariant derived category for the $\g$-conjugation action of $G$ on itself.

Now suppose $\g_1,\g_2:G\to G$ are automorphisms of $G$ that commute with each other. Then the following diagram commutes:

\beq\xymatrixcolsep{7pc}\xymatrix{
G\times G\ar[r]^-{\g_1\operatorname{-conjugation}}\ar[d]_{\g_2\times \g_2} & G\ar[d]^{\g_2}\\
G\times G\ar[r]^-{\g_1\operatorname{-conjugation}} & G.\\
}\eeq
Hence pullback by $\g_2$ induces a functor $\g_2^*:\D_G^{\g_1}(G)\to \DgaG$. We let $\DGab$ denote the equivariantization. Namely the objects of $\DGab$ are pairs $(M,\psi)$ where $M\in \DgaG$ and $\psi:\g_2^*(M)\xto{\cong} M$ is an isomorphism in $\DgaG$.

Let us define 
\beq 
\Rab:=\{(g,h)\in G\times G|h\g_2(g)\g_1(h)^{-1}=g\}.
\eeq
We equip $\Rab$ with an action of $G$ by setting ${ }^x(g,h)=(xg\g_1(x)^{-1}, xh\g_2(x)^{-1})$ for $x\in G$ and $(g,h)\in \Rab$. We let $\Fun_G(\Rab)$ denote the space of $\Qlcl$-valued $G$-invariant functions on $\Rab$. In our applications in this paper, one or both of the automorphisms will be Frobenius maps.

\brk\label{gammaequiv}
For $(g,h)\in \Rab$, we can check that ${ }^{g^{-1}}(g,h)=(\g_1(g),\g_1(h))$ and ${ }^{h^{-1}}(g,h)=(\g_2(g),\g_2(h))$.
\erk

\bex\label{gamfrob}
Let $G$ be connected, let $F:G\to G$ be a Frobenius map and let $\g:G\to G$ be an automorphism that commutes with $F$. Then the quotient set $G\backslash R_{\g,F}$ is finite and in fact can be naturally identified with the set $G^{F}/\sim_{\g}$. Note that in this case, $\g$ induces an automorphism of the finite group $G^{F}$. To construct this bijection, we note that since $G$ is connected, the $G$-orbit of each $(g,h)\in R_{\g,F}$ contains an  element of the form $(x,1)$ by Lang's Theorem. Also we see that $(x,1)\in R_{\g,F}$ if and only if $x\in G^{F}$. Moreover the $\g$-twisted conjugacy class of $x$ in $G^{F}$ is uniquely determined by the $G$-orbit of $(g,h)$. This defines the desired natural bijection. Hence in this case we have a natural identification of the spaces $\Fun_G(R_{\g,F})\cong \Fun(G^{F}/\sim_{\g})$. Some cases of particular interest to us are when $\g=\id_G$ or $\g$ is another Frobenius.
\eex

\brk 
In this paper we wish to work with possibly disconnected groups. In the general situation it is more natural to work with the space $\Fun_G(R_{\g,F})$ rather than the space $\Fun(G^F/\sim_\g)$ since the former takes into account all the pure inner forms.
\erk

\bex\label{fungandfunr}
Let $G$ be any algebraic group equipped with a Frobenius $F$. In this case the $G$ orbits in  $R_{\id_G,F}$ are in a natural bijection with the sets conjugacy classes in all the pure inner forms $G^{gF}$. Hence we have a natural identification $\FunG=\Fun_G(R_{\id_G,F})$. We will use these two notations interchangeably.
\eex

\subsection{Some auxiliary maps}\label{s:auxmaps}

The following simple lemma is key in defining the Shintani descent maps for general algebraic groups.
\blem\label{keylemma}
(i) Swapping the factors defines a $G$-equivariant isomorphism 
\beq
\tau:\Rab \to \Rba.
\eeq
(ii) We have a $G$-equivariant twisting isomorphisms 
\beq\label{t1def}
t_1:\Rab \to R_{\g_1\g_2, \g_2}
\eeq
defined by $t_1(g,h)=(h\g_2(g), h)=(g\g_1(h),h)$ and
\beq\label{t2def}
t_2:\Rab \to R_{\g_1,\g_1\g_2}
\eeq
defined by $t_2(g,h)=(g, h\g_2(g))=(g,g\g_1(h))$.
\elem

\bex\label{comfrob}
Suppose that $G$ is connected and that $F_1,F_2:G\to G$ are commuting Frobenius maps. By Example \ref{gamfrob} $G\backslash R_{F_1,F_2}=G^{F_2}/\sim_{F_1}$ and $G\backslash R_{F_2,F_1}=G^{F_1}/\sim_{F_2}$ and the swapping map $\tau$ induces an identification $G^{F_2}/\sim_{F_1}\cong G^{F_1}/\sim_{F_2}$ which agrees with the identification defined in \cite[\S1.1]{S} up to an inversion. This identification was used in \cite{S} to define Shintani descent for connected groups. 
\eex

\bex
Let $G$ be connected and let $F:G\to G$ be a Frobenius. Then by Example \ref{gamfrob} we have $G\backslash R_{\id_G,F}=G^F/\sim$. In this case the $G$-isomorphism $t_2:R_{\id_G,F}\to \Rif$ induces the twisting map $N_1:(G^F/\sim)\rar{} (G^F/\sim$) that is studied in \cite[\S4]{S} and that we introduced in (\ref{twistcon}). This twisting map is a map which permutes the rational conjugacy classes of $G^F$ within each geometric conjugacy class.
\eex

\subsection{The sheaf-function correspondence}\label{s:sfc}

Given an object $(M,\psi)\in \DGab$ we will define the associated ``trace function'' $T_{M,\psi}\in \Fun_G(\Rab)$. For any $g,h\in G$, the equivariance structure $\phi_M$ of $M\in \DgaG$ and the $\g_2$-equivariance structure $\psi$ give us isomorphisms of stalks 
\beq
M_{h\g_2(g)\g_1(h)^{-1}}\xto{\phi_M(h^{-1}, h\g_2(g)\g_1(h)^{-1})} M_{\g_2(g)}\xto{\psi(g)}M_{g}.
\eeq
 If $(g,h)\in \Rab$, then this composition is an automorphism of the stalk $M_g\in D^b(\Vec)$. Hence we can define $T_{M,\psi}(g,h)$ to be the trace of this composition.

\bex\label{trfun}
Let $G$ be any algebraic group. Set $\g_1=\id_G$ and $\g_2$ to be an $\F_q$-Frobenius map $F:G\to G$. In this case, the construction of the trace functions defined here reduces to the notion of the ``trace function'' associated with an $F$-equivariant object of $\DG$ that was studied in \cite[\S2.4.8]{De}. Recall from Example \ref{fungandfunr} that in this case the space $\Fun_G(R_{\id_G,F})$ can be identified with the space $\FunG$ of class functions on all pure inner forms of $G^F$. We refer to {\it loc. cit.} for details.
\eex

\bex\label{charfun}
Let $G$ be any algebraic group, but now let $\g_1$ be a Frobenius $F:G\to G$ and let $\g_2=\id_G$.  In this case, our construction takes an object in $\D^F_G(G)^{\id_G}$ and produces a function in $\Fun_G(R_{F,\id_G})$. Let $(M,\id_M)\in \D^F_G(G)^{\id_G}$. This gives us a function $\tau^*T_{M,\id_M}\in \Fun_G(R_{\id_G,F})$ where $\tau:R_{\id_G,F}\to R_{F,\id_G}$ is the swapping map. Then it is easy to check that $\tau^*T_{M,\id_M}=\bar{\chi_M}$, where $\chi_M$ is the character of $M\in \D^F_G(G)$ as defined in \cite[\S2.4.7]{De}. 
\eex

\subsection{Algebraic groups over $\Fq$ equipped with an $\Fq$-automorphism}\label{s:agof}
In this section we study algebraic groups defined over $\Fq$ equipped with an automorphism defined over $\Fq$. Let $G$ be any algebraic group equipped with an $\Fq$-Frobenius $F:G\to G$ and an automorphism $\g:G\to G$ commuting with $F$. Let us now understand the map of $\k$-schemes $R_{\g,F}\rar{p_2}G$. For this it is convenient to identify $G$ with the coset $G\g\subset G\rtimes \g^\Z$. Then the conjugation action of $G$ on $G\g$ corresponds to the $\g$-twisted conjugation on $G$. For each $h\in G$, we have the Frobenius map $hF:G\g\to G\g.$ Under the identification with $G$, this  modified Frobenius maps $g\mapsto hF(g)\g(h)^{-1}$. Note that this modified Frobenius does not respect the group structure of $G$. In other words, we have
\blem\label{l:silly}
In the notation above, $(g,h)\in R_{\g,F}$ if and only if $g\in (G\g)^{hF}$.
\elem
\brk\label{r:silly}
Hence we can identify the spaces $\Fun_G(R_{\g,F})$ and $\Fun([G\g],F)$, the space of class functions on the $\Fq$-points of all pure inner forms of the coset $G\g$.
\erk

Let us now study the category $\D^F_G(G)^\g$ and the space $\Fun_G(R_{F,\g})$. There are only finitely many $G$-orbits in $R_{F,\g}$ and hence $\Fun_G(R_{F,\g})$ is finite dimensional. For $(g,h)\in R_{F,\g}$ its stabilizer, 
\beq
\Stab_G(g,h)=\{x\in G^{gF}|h\g(x)h^{-1}=x\}
\eeq
is finite. We define an Hermitian inner product on $\Fun_G(R_{F,\g})$ as follows:
\beq
\<f_1,f_2\>:=\sum\limits_{\<(g,h)\>\in G\backslash R_{F,\g}}{\frac{f_1(g,h)\bar{f_2(g,h)}}{|\Stab_G(g,h)|}}.
\eeq 

As stated earlier, we need to take into consideration all pure inner forms of $F$. For each $g\in G$ we have the following commutative diagram
\beq\xymatrix{
G\ar[r]^-{\g}\ar[d]_{\ad(g)\circ F} & G\ar[d]^{\ad(\g(g))\circ F}\\
G\ar[r]^-{\g} & G.\\
}\eeq
Hence $\g$ induces an isomorphism $\g:G^{gF}\xto{\cong} G^{{\g(g)}F}$ and hence it induces a bijection 
\beq
\g^*:\Irrep(G^{{\g(g)}F})\xto{\cong}\Irrep(G^{gF}).
\eeq 
We recall from \cite{De} that the set $\Irrep(G^{gF})$ is canonically determined by the $F$-conjugacy class of $g$ in $H^1(F,G)$. Hence we see that $\g$ induces a bijection 
\beq
\g^*:\Irrep(G,F)\to \Irrep(G,F).
\eeq

\subsection{Irreducible representations fixed by automorphisms}\label{s:irrepgam}
In this section we study the fixed point set $\Irrep(G,F)^\g$ of the bijection $\g^*:\Irrep(G,F)\to \Irrep(G,F)$ defined above. Let us recall a simple observation from \cite[\S2.4.2]{De}. (See \S\ref{s:tconjact} for the notation.)
\blem\label{dfg}
We have an equivalence of triangulated categories
$$\D^F_G(G)\cong \bigoplus\limits_{\<g\>\in H^1(F,G)}D^b\Rep(G^{gF}).$$
The category $\D^F_G(G)$ is in fact a semisimple abelian category whose (isomorphism classes of) simple objects are parameterized by $\Irrep(G,F)\times \Z$, where the integers keep track of the degree shift. If $W\in \Irrep(G^{gF})\subset \Irrep(G,F)$ then we have the associated local system $W_{loc}\in \D^F_G(G)$ supported on the $F$-twisted conjugacy class of $g$.
\elem

We note that $\g$ acts on the set $H^1(F,G)$ which parameterizes the pure inner forms and only the irreducible representations of $\g$-fixed pure inner forms can lie in $\Irrep(G,F)^\g$. Hence suppose that $g\in G$ represents a $\g$-fixed $F$-twisted conjugacy class i.e. there exists an $h\in G$ such that $g=h\g(g)F(h)^{-1}$, or in other words an $h$ such that $(g,h)\in R_{F,\g}$. We have an automorphism 
\beq
h\g:G^{gF}\rar{\g} G^{{\g(g)}F}\xto{\ad(h)}G^{gF}.
\eeq
We note that a different choice of $h\in G$ only changes the automorphism above by an inner automorphism of $G^{gF}$. Hence the action of $h\g$ on $\Irrep(G^{gF})$ is independent of the choice of $h$. Let $W\in \Irrep(G^{gF})\subset \Irrep(G,F)$ and consider the local system $W_{loc}\in \D_G^F(G)$ supported on the $F$-twisted conjugacy class of $g$. Then $W\in \Irrep(G,F)^\g$ if and only if there exists an intertwiner $\psi:(h\g)^*W\rar{\cong}W$. Then we can define an object $(W_{loc},\psi)\in \D^F_G(G)^\g$. Taking the associated trace function, we get $T_{W,\psi}\in \Fun(G^{gF}/\sim_{h\g})\subset\Fun_G(R_{F,\g})$ (see Lemma \ref{onbasis}). In fact we can choose the intertwiner $\psi$ in such a way that we obtain an irreducible representation of the finite group $G^{gF}\<h\g\>$, where $\<h\g\>$ is the finite cyclic group generated by the automorphism $h\g$ of $G^{gF}$. For each $W\in \Irrep(G,F)^\g$ we fix such an intertwiner $\psi_W$ as above. In this case we have $|T_{W,\psi_W}|=1$ and $T_{W,\psi_W}$ is well-defined by $W$ up to scaling by a root of unity.

In other words we have:
\blem\label{onbasis}
Let $W\in \Irrep(G,F)$. Then $W\in \Irrep(G,F)^\g$ if and only if $\g^*W_{loc}\cong W_{loc}$. If such a $W$ lies in $\Irrep(G^{{gF}})^{h\g}$, we can choose an intertwiner $\psi_W:(h\g)^*W\rar{\cong}W$ such that we obtain a representation of $G^{gF}\<h\g\>$ on $W$ extending the original representation of $G^{gF}$. Such a $\psi_W$ is well-defined up to scaling by roots of unity. Further we have an orthogonal decomposition 
\beq
\Fun_G(R_{F,\g})=\bigoplus\limits_{\<g\>\in H^1(F,G)|\g\<g\>=\<g\>}\Fun(G^{gF}/\sim_{h\g}).
\eeq 
The set 
\beq
\{T_{W,\psi_W}|W\in \Irrep(G^{gF})^{h\g}\}\subset \Fun(G^{gF}/\sim_{h\g})
\eeq is an orthonormal basis and hence the set 
\beq
\{T_{W,\psi_W}|W\in \Irrep(G,F)^{\g}\}\subset \Fun_G(R_{F,\g})
\eeq is an orthonormal basis.
\elem

\subsection{Irreducible representations fixed by Frobenius}\label{s:irrepfrob}

Let us now refine the results of the previous section in the special case which is of interest from the point of view of Shintani descent. Let $m$ be any positive integer and let us consider the set $\Irrep(G,F^m)^F$ which appears on one side of the Shintani descent map. If $W\in \Irrep(G,F^m)^F$ we will place some further restrictions on the choice of the intertwiner $\psi_W$ such that with these restrictions the different choices differ by $m$-th roots of unity. Hence with each $W\in \Irrep(G,F^m)^F$ we will be able to associate a function $T_{W,\psi_W}\in \Fun_G(R_{F^m,F})$ which will be unique up to scaling by $m$-th roots of unity.

Suppose $W\in \Irrep(G,F^m)^F$ where $W$ is an irreducible representation of some pure inner form $G^{gF^m}$. By \S\ref{s:irrepgam} there must be an $h\in G$ such that $(g,h)\in R_{F^m,F}$. We can check that $hF(h)\cdots F^{m-1}(h)g^{-1}\in G^{gF^m}\cap G^{hF}$. Consider the group automorphism $hF:G^{gF^m}\to G^{gF^m}$ from \S\ref{s:irrepgam}. Then $W\in \Irrep(G^{gF^m})^{hF}$. Let $\t{G^{gF^m}\<hF\>}$ be the quotient of the semidirect product $G^{gF^m}\rtimes (hF)^{\Z}$ modulo the relation $G^{gF^m}\ni hF(h)\cdots F^{m-1}(h)g^{-1}=(hF)^m$. This gives us an extension
\beq\label{extbycyc}
0\to G^{gF^m}\to  \t{G^{gF^m}\<hF\>} \to \Z/m\Z \to 0.
\eeq
Let us note that the semidirect product $G^{gF^m}\rtimes (hF)^{\Z}$ as well as the extension above do not depend (up to a canonical isomorphism) on the choice of $h$ since the possible choices all lie in the coset $G^{gF^m}h$. For example, a different choice of $h$ merely amounts to a different splitting of the semidirect product.

Since $W\in \Irrep(G^{gF^m})^{hF}$ we can extend the action to obtain a representation of $\t{G^{gF^m}\<hF\>}$ on $W$. This extension amounts to choosing a suitably normalized intertwiner $\psi_W: (hF)^*W\to W$ of $G^{gF^m}$-representations. We thus obtain an object $(W_{loc},\psi_W)\in \D^{F^m}_G(G)$ and the corresponding trace function $T_{W,\psi_W}\in \Fun_G(R_{F^m,F})$. It is clear that the different choices of extensions of the representation as above, namely the different choices of normalization of the intertwiner $\psi_W$ differ by scaling by $m$-th roots of unity. (See also Remark \ref{r:altpsi} for another equivalent characterization of $\psi_W$.) Hence we have proved that given $W\in  \Irrep(G,F^m)^F$, we can obtain the function $T_{W,\psi_W}\in \Fun_G(R_{F^m,F})$ as above, uniquely determined up to scaling by $m$-th roots of unity. By Lemma \ref{onbasis}, the set $\{T_{W,\psi_W}\}_{W\in \Irrep(G, F^m)^{F}}$ is an orthonormal basis of $\Fun_G(R_{F^m,F})$.

\subsection{Convolution of functions and module structure}\label{s:convfun}
Let $F:G\to G$ be a Frobenius map and let $\g_1,\g_2:G\rar{} G$ be group automorphisms that commute with the Frobenius. In this section we define convolution products (which are associative)
\beq\label{convfun}
\Fun_G(R_{\g_1,F})\times \Fun_G(R_{\g_2,F}) \rar{} \Fun_G(R_{\g_1\g_2,F}).
\eeq

If $(g_i,h)\in R_{\g_i,F}$ for $i=1,2$, then we can check that $(g_1\g_1(g_2),h)\in R_{\g_1\g_2,F}$. Then for $f_i\in \Fun_G(R_{\g_i,F})$ we define
\beq
f_1\ast f_2(g,h)=\sum\limits_{\substack{{g_1\g_1(g_2)=g}\\{(g_i,h)\in R_{\g_i,F}}}}f_1(g_1,h)f_2(g_2,h) 
\eeq
for each $(g,h)\in R_{\g_1\g_2,F}$. Note that the sum is finite (for example using Lemma \ref{l:silly}) and hence the convolution is well defined. If is easy to check that $f_1\ast f_2\in \Fun_G(R_{\g_1\g_2,F})$ and also that the convolution is associative. 
\brk
This convolution of functions may perhaps be more transparent using the identification of Remark \ref{r:silly}.
\erk
\brk\label{r:module}
Using the convolution, each function space $\Fun_G(R_{\g,F})=\Fun([G\g],F)$ in fact becomes a module over the commutative algebra $\FunG$.
\erk

\subsection{Convolution of complexes and sheaf-function correspondence}\label{s:convandsf}
Let $\g_1,\g_2,F$ be as in the previous section. Then using the ``coset multiplication'' 
\beq
\mu_{\g_1,\g_2}:G\times G\to G
\mbox{ defined by $(g_1,g_2)\mapsto g_1\g_1(g_2)$}
\eeq we can define a convolution (with compact supports) 
\beq
\D^{\g_1}_G(G)\times \D^{\g_2}_G(G)\rar{}\D_G^{\g_1\g_2}(G).
\eeq It may be useful to think of this as the standard convolution $\D_G(G\g_1)\times \D_G(G\g_2)\rar{}\D_G(G\g_1\g_2)$. The map $\g_1:G\rar{}G$ induces the equivalence 
\beq
\g_1:={\g_1^{-1}}^*={\g_1^*}^{-1}:\D_G^{\g_1}(G)\rar{\cong}\D_G^{\g_1\g_2\g_1^{-1}}(G).
\eeq 
Moreover for $M_i\in\D^{\g_i}_G(G)$ we have functorial crossed braiding isomorphisms
\beq
\beta_{M_1,M_2}:M_1\ast M_2\rar{\cong}\g_1(M_2)\ast M_1.
\eeq
We also have the induced convolution
\beq
\D^{\g_1}_G(G)^F\times \D^{\g_2}_G(G)^F\rar{}\D_G^{\g_1\g_2}(G)^F.
\eeq

Now the sheaf-function correspondence is compatible with pullbacks and pushforwards with compact supports (see also \cite[\S2.2]{De}). Using this, we obtain:
\bprop\label{p:convandsf}
The convolution of complexes (as defined above) is compatible with the convolution of functions (as defined in \S\ref{s:convfun})  under the sheaf-function correspondence. Namely, let $M_i\in \D^{\g_i}_G(G)^F$ and let $T_{M_i}\in \Fun_G(R_{\g_i,F})$ be the associated trace function.  Then we have 
\beq
T_{M_1}\ast T_{M_2}=T_{M_1\ast M_2}.
\eeq
\eprop

\section{Shintani descent for general algebraic groups}\label{s:sdgag}
In this section we define Shintani descent using the constructions from \S\ref{s:prelcon}.

\subsection{Definition of Shintani descent}\label{s:defsd}
In this section we apply our previous results and constructions and define Shintani descent in general. Let $G$ be any algebraic group over $\k$. Let $F:G\to G$ be an $\Fq$-Frobenius map. Let $m$ be a positive integer. Note that we have a $G$-equivariant isomorphism (which we will call the inverse norm map)
\beq
N_m^{-1}:R_{1,F}\rar{\cong}R_{F^m,F}
\eeq 
given by $R_{1,F}\ni (g,h)\mapsto (ghF(h)\cdots F^{m-1}(h),h)\in R_{F^m,F}$, i.e. it is the composition of the twists (\ref{t1def})
\beq
R_{1,F}\xto{t_1}R_{F,F}\xto{t_1} R_{F^2,F}\xto{t_1}\cdots R_{F^m,F}.
\eeq
Hence we obtain an isomorphism of Hermitian inner product spaces 
\beq
{N_m^{-1}}^*:\Fun_G(R_{F^m,F})\rar{\cong}\Fun_G(R_{1,F})=\FunG.
\eeq 

\blem\label{l:modiso}
For each positive integer $m$, the map $t_1^*:\Fun_G(R_{F^m,F})\rar{}{\Fun_G(R_{F^{m-1},F})}$ is an isomorphism of $\FunG$-modules (see also Remark \ref{r:module}). Hence the map ${N_m^{-1}}^*$ defined above is a $\FunG$-module isomorphism.
\elem
\bpf
It is enough to check that $t_1^*$ preserves the module structure. Moreover, it is easy to check that $f_1\ast t_1^*f_2=t_1^*(f_1\ast f_2)$ for $f_1\in \FunG$ and $f_2\in \Fun_G(R_{F^m,F})$.
\epf

Given $W$ in $\Irrep(G,F^m)^F$, we have the associated trace function $T_{W,\psi_W}\in \Fun_G(R_{F^m,F})$ as in \S\ref{s:irrepfrob}.
\bdefn
We define Shintani descent 
\beq
\Sh_m:\Irrep(G,F^m)^F\hookrightarrow\FunG=\Fun_G(R_{\id_G,F}) 
\eeq
by $\Sh_m(W)={N_m^{-1}}^*(T_{W,\psi_W})$ for each $W\in \Irrep(G,F^m)^F$. By Lemma \ref{onbasis}, the image of Shintani descent $\{\Sh_m(W)\}_{W\in \Irrep(G, F^m)^{F}}\subset \FunG$ is an orthonormal basis known as the $m$-th Shintani basis.
\edefn
\brk
By \S\ref{s:irrepfrob}, for each $W\in \Irrep(G,F^m)^F$, $\Sh_m(W)$ is uniquely determined up to scaling by an $m$-th root of unity.
\erk

\subsection{The twisting operator}\label{s:twistoper}
In this section we study a special case of Shintani descent when $m=1$. In this case Shintani descent is a map
\beq
\Sh_1:\Irrep(G,F)=\Irrep(G,F)^F\hookrightarrow\FunG.
\eeq
However as we will see now, this is not in general the natural inclusion obtained by taking characters, rather it is a twist of the natural inclusion.

Consider any object $M\in \D_G^F(G)$. Then as defined in \cite[\S2.4.7]{De} (see also Example \ref{charfun}) we have its character $\chi_M\in \Fun_G(R_{\id,F})$. Let $\phi_M$ denote the $F$-twisted conjugation equivariance structure associated with $M$. We recall that for $(g,h)\in R_{\id,F}$, $\chi_M(g,h):=\tr(M_h\xto{\phi_M(g,h)}M_h)$. We further set 
\beq
\psi_M(g):=\phi_M(g,F(g)):M_{F(g)}\rar{\cong}M_g
\eeq and obtain an object $(M,\psi_M)\in \D_G^F(G)^F$. Let us compute its associated trace function $T_{M,\psi_M}\in \Fun_G(R_{F,F})$. Suppose $(g,h)\in R_{F,F}$ i.e. $hF(g)F(h)^{-1}=g$. Then $T_{M,\psi_M}(g,h)$ is the trace of the automorphism
\beq
M_g\xto{\phi_M(h^{-1}, g)} M_{F(g)}\xto{\phi_M(g,F(g))}M_{g}.
\eeq
Hence $T_{M,\psi_M}(g,h)=\tr(\phi_M(gh^{-1},g))=\chi_M(gh^{-1},g)$. Now let $(g,h)\in R_{\id,F}$. Then using Lemma \ref{keylemma} we see that
\beq
{N_1^{-1}}^*T_{M,\psi_M}(g,h)=T_{M,\psi_M}(gh,h)=\chi_M(g,gh)=t_2^*\chi_M(g,h),
\eeq where $t_2$ is as in (\ref{t2def}). Hence we see that if $W\in \Irrep(G,F)$, then we can canonically define $\Sh_1(W)$ and we have $\Sh_1(W)=t_2^*\chi_W$ where $\chi_W\in \Fun_G(R_{\id,F})$ denotes the character of $W$.

\subsection{Sheaf-function correspondence and twists}\label{s:twistandsf}
As we have seen before, the twist $t_2:R_{\id,F}\rar{}R_{\id,F}$ from (\ref{t2def}) induces a unitary operator
\beq
t_2^*:\FunG\rar{}\FunG.
\eeq

In general an irreducible character $\chi_W$ for $W\in \Irrep(G,F)$ may not be an eigenvector for the twisting operator defined above. However as we will see below, if $(C,\psi)\in \D_G(G)^F$ is a conjugation equivariant Weil complex such that $\End(C)=\Qlcl$, then the trace of Frobenius function $T_{C,\psi}\in \FunG$ is an eigenvector for the twisting operator.

First we recall from \cite{BD} that the category $\DG$ is equipped with a twist $\theta$ i.e. an automorphism of the identity functor satisfying some properties. Hence for each $C\in \DG$ we have the twist $\theta_C:C\rar{\cong}C$. We now prove
\bprop\label{p:sftwist}
Let $(C,\psi)\in \DG^F$ be a Weil complex such that $\theta_C$ is a scalar. Then \beq t_2^*T_{C,\psi}=\theta_C^{-1}\cdot T_{C,\psi}.\eeq
\eprop
\bpf
Let $(g,h)\in R_{\id,F}$. Then we want to prove that
\beq\label{eqsftwist}
T_{C,\psi}(g,gh)=\theta_C^{-1} T_{C,\psi}(g,h).
\eeq 
We have a commutative diagram
\beq\xymatrixcolsep{5pc}\xymatrix{
M_g\ar[r]^-{\phi_C(h^{-1},g)}\ar[d]_{\phi_C(g,g)=\theta_C} & M_{F(g)}\ar@{=}[d]\ar[r]^-{\psi(g)} & M_g\ar@{=}[d]\\
M_g\ar[r]^-{\phi_C(h^{-1}g^{-1},g)} & M_{F(g)}\ar[r]^-{\psi(g)} & M_g,\\
}\eeq
where $\phi_C$ denotes the $G$-equivariance structure associated with $C\in \DG$. The equality $\phi_C(g,g)=\theta_C$ follows from the definition of the twist $\theta_C$ (see \cite{BD}). The trace of the top row is $T_{C,\psi}(g,h)$ and the trace of the bottom row is $T_{C,\psi}(g,gh)$. Hence (\ref{eqsftwist}) holds.
\epf

\section{Traces and twists in module categories}\label{s:twisttrace}
We have seen that the category $\D^F_G(G)$ is a $\DG$-module category. In this section, we will study some additional structures on $\D^F_G(G)$ known as twists and traces.

\subsection{Twists in the category $\D^F_G(G)$}\label{modcatandtwist}
We know that the category $\D^F_G(G)$ is a $\DG$-module category and that the category $\DG$ is equipped with a twist $\theta:\id_{\DG}\to \id_{\DG}$. Explicitly, if $C\in \DG$ and $\phi_C$ is the conjugation equivariance structure associated with $C$, then $\theta_C(g)=\phi_C(g,g):C_g\rar\cong C_g$. We will now construct an analgous twist for the module category $\D^F_G(G)$, namely a natural transformation $\theta^F:\id_{\D^F_G(G)}\to F|_{\D^F_G(G)}={F^{-1}}^*|_{\D^F_G(G)}$.

Let $M\in \D^F_G(G)$ and let $\phi_M$ denote the $G$-equivariance structure of $M$. Then $\theta^F_{M}:M\to F(M)$ is defined by 
\beq
\theta^F_{M}(g)=\phi_M(F^{-1}(g),g):M_g\rar\cong M_{F^{-1}(g)}=F(M)_g.
\eeq

It is clear that the Frobenius preserves the twists, namely
\beq
F(\theta_C)=\theta_{F(C)} \mbox{ for C $\in \DG$ and }
\eeq
\beq
F(\theta^F_M)=\theta^F_{F(M)} \mbox{ for M $\in \D^F_G(G)$. }
\eeq

\brk\label{r:altpsi}
For each positive integer $m$, we have a twist $\theta^{F^m}$ in the $\DG$-module category $\D^{F^m}_G(G)$. If we have $W\in \Irrep(G,F^m)^F$, with the choice of $\psi_W$ as in \S\ref{s:irrepfrob}, then the composition
\beq
{F^m}^*W_{loc}\xto{{F^{(m-1)}}^*(\psi_W)}{F^{(m-1)}}^*W_{loc}\to\cdots F^*W_{loc}\xto{\psi_W}W_{loc}
\eeq equals ${F^m}^*\theta^{F^m}_{W_{loc}}$
\erk

It is also straightforward to check that the twist $\theta^F$ is compatible with the twist in the category $\DG$:
\blem\label{modcattwist}
Let $C\in \DG$ and let $M\in \D^F_G(G)$. Then $\theta^F_{C\ast M}$ is equal to the composition
\beq
C\ast M\xto{\beta_{C,M}} M\ast C\xto{\beta_{M,C}} F(C)\ast M \xto{\theta_{F(C)}\ast \theta^F_M}F(C)\ast F(M)\rar{} F(C\ast M).
\eeq
Moreover, all the twists $\theta^{F^a}$ in $\D^{F^a}_G(G)$ for $a\in \Z$ are compatible with each other, namely if $L\in \D^{F^a}_G(G)$, $M\in \D^{F^b}_G(G)$ then $\theta^{F^{a+b}}_{L\ast M}$ equals the composition
\beq
L\ast M\xto{\beta_{L,M}} F^a(M)\ast L\xto{\beta_{F^a(M),L}} F^b(L)\ast F^a(M) \xto{\theta^{F^a}_{F^b(L)}\ast \theta^{F^b}_{F^a(M)}}F^{a+b}(L)\ast F^{a+b}(M)\rar{} F^{a+b}(L\ast M).
\eeq
\elem

%For each $m\in \Z$, we have a weak duality functor $\f{D}^{-}:\D^{F^m}_G(G)\rar{\cong} \D^{F^{-m}}_G(G)^{\h{op}}$ (cf. \cite{BD}) which commutes with the Frobenius. Furthermore, for $M\in \D^{F^m}_G(G)$ we have 
%\beq
%\f{D}^-\theta^{F^m}_M=\theta^{F^{-m}}_{F^m(\f{D}^-M)}:F^m(\f{D}^-M)\cong \f{D}^-F^m(M)\rar{\cong}\f{D}^-M.
%\eeq

\subsection{Compatibility of twists and traces}
In \cite[\S2.4.5]{De}, we have defined a normalized trace $\tr_F$ in the category $\D^F_G(G)$ which assigns a number in $\Qlcl$ to each endomorphism in $\D^F_G(G)$. We now prove:
\bthm\label{t:twisttrace}
For each object $M\in \D^F_G(G)$ we define the automorphism $\nu_{\theta^F_M}$ in $\D^{F^2}_G(G)$ to be the composition
\beq
\nu_{\theta^F_M}:M\ast M\xto{\beta_{M,F^*M}^{-1}} M\ast F^*M\xto{\id_M\ast F^*(\theta^F_M)}M\ast M.
\eeq
Then we have
\beq
\tr_F(\id_M)=\tr_{F^2}(\nu_{\theta^F_M}).
\eeq
\ethm
\bpf
As usual let $\phi_M$ denote the $G$-equivariance structure associated with $M\in \D^F_G(G)$. We will use the integral symbol ($\int$) to denote pushforward with compact supports; if $\E\in \D(X)$, we use $\int\limits_{X}\E$ to denote $R\Gamma_c(X,\E)$. Let $t\in G$. Using the definition of the convolution, the crossed braidings and the twists, we see that the stalk
$$\nu_{\theta^F_m}(t):(M\ast M)_t\rar{} (M\ast M)_t$$ is given by the composition
$$(M\ast M)_t=\int\limits_{h_1F(h_2)=t} M_{h_1}\otimes M_{h_2} $$
$$\xto{{\int}\phi_M(h_2^{-1},h_1)\otimes \id_{M}}\int\limits_{h_1F(h_2)=t} M_{h_2^{-1}h_1F(h_2)}\otimes M_{h_2}$$
$$=\int\limits_{h_1F(h_2)=t} F^*M_{F^{-1}(h_2^{-1}h_1)h_2}\otimes M_{h_2}=(M\ast F^*M)_t$$
$$\xto{F^*(\theta_M^F)\ast\id_M=\int\phi_M(F^{-1}(h_2^{-1}h_1)h_2,h_2^{-1}h_1F(h_2))\otimes \id_M}\int\limits_{h_1F(h_2)=t} M_{F^{-1}(h_2^{-1}h_1)h_2}\otimes M_{h_2}$$
%$$\xto{{\displaystyle\int}\left(\psi_L(F^{-1}(F^{-m}(h_2)^{-1}h_1h_2))\circ\phi_L(F^{-m}(h_2)^{-1},h_1)\right)\otimes \left(\psi_M(F^{-1}(h_2^{-1}h_1^{-1}F^{-m}(h_2))h_1h_2)\circ\phi_M(F^{-1}(F^m(h_2^{-1}h_1^{-1})h_2), h_2)\right)}$$
$$\xto{\tau=\h{swap}}\int\limits_{h_1F(h_2)=t}M_{h_2} \otimes M_{F^{-1}(h_2^{-1}h_1)h_2} $$
$$=(M\ast M)_t.$$

Consider the automorphism of the antidiagonal $\bar{\Delta}_t:=\{(h_1,h_2)\in G\times G|h_1F(h_2)=t\}$ that appears implicitly in the stalk of the automorphism $\zeta_{L,M}$ as above, namely
$$(h_1,h_2)\mapsto (h_2,F^{-1}(h_2^{-1}h_1)h_2).$$
If we identify $\bar{\Delta}_t$ and $G$ using the second projection, then the inverse of this automorphism corresponds to the Frobenius automorphism of $G$ (as a scheme) defined by $h_2\mapsto tF(h_2^{-1})$. The fixed point set of this Frobenius map is  the set $\{(h,h)|hF(h)=t\}\subset\bar{\Delta}_t$. Hence by the Grothendieck-Lefschetz trace formula (trace of the induced map on cohomology equals the sum of traces of the stalks of the map over all fixed points, see also \cite[Lemma 4.4(iii)]{B}), we deduce that 
$$
\tr(\nu_{\theta^F_M}(t))=\sum\limits_{hF(h)=t}\tr(\tau\circ(\phi_M(1,h)\otimes \id):M_h\otimes M_h\xto{\cong}M_h\otimes M_h)
$$
$$
=\sum\limits_{hF(h)=t}\tr(\tau:M_h\otimes M_h\xto{\cong}M_h\otimes M_h)
$$
$$
=\sum\limits_{hF(h)=t}\dim M_h.
$$
Now $M\ast M$ is an object of $\D_G^{F^{2}}(G)\cong \bigoplus\limits_{\substack{<t> \in \\ H^1(F^{2},G)}}D^b\h{Rep}(G^{tF^2})$. Hence by definition of the trace $\tr_{F^{2}}$ on $\D^{F^{2}}_G(G)$ (see \cite[\S2.4.5]{De}) we have 
\beq
\tr_{F^{2}}(\nu_{\theta^F_M})=\sum\limits_{\substack{<t> \in \\ H^1(F^{2},G)}}\frac{\tr(\nu_{\theta^F_M}(t))}{|G^{tF^2}|}
\eeq
\beq
=\sum\limits_{\substack{<t> \in \\ H^1(F^{2},G)}}\sum\limits_{hF(h)=t}\frac{\dim M_h}{|G^{tF^2}|}.
\eeq
Next note that $\Delta:=\{(h,h)|h\in G\}$ is in fact a $G$-invariant subscheme of $R_{F,F}$ and $(h,h)\mapsto (h,hF(h))$ under  the $G$-equivariant isomorphism $R_{F,F}\rar{t_2}R_{F,F^2}$. Hence we obtain
\beq
\tr_{F^2}(\nu_{\theta^F_M}(t))=\sum\limits_{\substack{<t> \in \\ H^1(F^{2},G)}}\sum\limits_{(h,t)\in t_2(\Delta)}\frac{\dim M_h}{|G^{tF^2}|}
\eeq
\beq
=\sum\limits_{\substack{<t> \in \\ H^1(F^{2},G)}}\sum\limits_{(h,t)\in t_2(\Delta)}\frac{\dim M_h}{|\Stab_G(h,t)|\cdot|G^{tF^2}\mbox{-orbit of }(h,t)|}
\eeq
\beq
=\sum\limits_{\substack{<(h,t)> \in G\backslash}t_2(\Delta)}\frac{\dim M_h}{|\Stab_G(h,t)|}
\eeq
\beq
=\sum\limits_{<(h,h)> \in G\backslash\Delta}\frac{\dim M_h}{|\Stab_G(h,h)|}
\eeq
\beq
=\sum\limits_{<h>\in H^1(G,F)}\frac{\tr(\id_{M_h})}{|G^{tF}|}
\eeq
\beq
\tr_F(\id_M)
\eeq
as desired.

\epf

\section{An inner product formula}\label{s:ipf}
For each positive integer $m$, we have the $m$-th Shintani basis $\{\Sh_m(W)\}_{W\in \Irrep(G, F^m)^{F}}$ of $\FunG$ defined in \S\ref{s:defsd}. On the other hand, we have the basis $\Irrep(G,F)$ consisting of the irreducible characters of all the pure inner forms. Our goal in this section is to describe the relationship between these two bases.

\subsection{Twists and sheaf-function correspondence}\label{s:crbr}
We know that the category $\D^F_G(G)$ is equivalent to the direct sum of the bounded derived categories of the representation categories of all pure inner forms. Let us now observe that we have a natural identification of all the categories $\D^F_G(G)^{F^m}$ for all integers $m\geq 0$.  
\blem\label{d:eta}
For each $m\geq 0$, we have an equivalence of triangulated categories
\beq
\D^F_G(G)^{\id_G}\rar{\cong} \D^F_G(G)^{F^m}
\eeq 
and hence the composition
\beq
\eta_m:\D^F_G(G)\rar{}\D^F_G(G)^{\id_G}\rar{\cong}  \D^F_G(G)^{F^m}
\eeq 
where the first functor is defined by $M\mapsto (M,\id_M)$.
\elem
\bpf
It is sufficient to construct an equivalence
\beq
\D^F_G(G)^{F^{m}}\rar{\cong} \D^F_G(G)^{F^{m+1}}
\eeq
for each $m\geq 0.$ Let $(M,\psi)\in \D^F_G(G)^{F^{m}}$, where $M\in \D^F_G(G)$ and $\psi:{F^m}^*M\rar{\cong}M$. On the other hand we have the twist $\theta^F_{{F^{(m+1)}}^*M}={F^{(m+1)}}^*\theta^F_M:{F^{(m+1)}}^*M\rar\cong {F^{m}}^*M$. Then $(M,\psi)\mapsto (M,\psi\circ{F^{(m+1)}}^*\theta^F_M)$ defines the desired equivalence.
\epf

\brk
From the proof, we see that if $M\in \D^F_G(G)$, then $\eta_m(M)=(M,\psi_{M,m})\in \D^F_G(G)^{F^m}$ where
\beq\label{e:eta}
\psi_{M,m}:{F^m}^*M\xto{{F^m}^*\theta^F_M}{F^{(m-1)}}^*M\xto{}\cdots F^*M\xto{F^*\theta^F_M}M.
\eeq
\erk

\blem\label{twistandsf}
The following diagram is commutative
\beq\xymatrixcolsep{2pc}\xymatrix{
\D^F_G(G)^{\id_G}\ar[r]^{\cong}\ar[d] & \D^F_G(G)^{F}\ar[d]\ar[r]^-{\cong} & \cdots\ar[r]^-{\cong} & \D^F_G(G)^{F^m}\ar[d]\ar[r]^-{\cong} &\\
\Fun_G(R_{F,\id_G})\ar[r]^-{\cong}\ar[d]^{\tau^*} & \Fun_G(R_{F,F})\ar[r]^-{\cong}\ar[d]^{\tau^*} & \cdots\ar[r]^-{\cong} &\Fun_G(R_{F,F^m})\ar[r]^-{\cong}\ar[d]^{\tau^*} & \\
\FunG\ar[r]^-{\cong} & \Fun_G(R_{F,F})\ar[r]^-{\cong} & \cdots\ar[r]^-{\cong} &\Fun_G(R_{F^m,F})\ar[r]^-{\cong} & ,\\
}\eeq
where the equivalences in the top row are as defined in the proof of the previous lemma, the isomorphisms in the bottom two rows are induces by the twists $t_2, t_1$ from (\ref{t2def}), (\ref{t1def})  respectively and the top vertical arrows are defined by the sheaf-function correspondence from \S\ref{s:sfc}.
\elem
\bpf
This is straighforward to check from the definitions.
\epf

\subsection{The inner product formula}\label{s:tipf}

Let $m\geq 0$ be an integer. Suppose that we have objects $L\in \D^{F^m}_G(G)$ and $M\in \D^F_G(G)$ then we have the following functorial crossed braiding isomorphisms (see \S\ref{s:convandsf})
\beq
\beta_{M,F^*L}^{-1}: L\ast M\rar{} M\ast F^*L
\eeq
and 
\beq
\beta_{L,{F^m}^*M}^{-1}: M\ast L\rar{} L\ast {F^m}^*M
\eeq
in the category $\D^{F^{m+1}}_G(G)$. If $L\in \D^{F^m}_G(G)^F$ and $M\in \D^F_G(G)^{F^m}$ then we obtain an automorphism
\beq\label{defzeta}
\zeta_{L,M}: L\ast M\xto{\beta_{M,F^*L}^{-1}} M\ast F^*L \xto{\beta_{F^*L,{F^m}^*M}^{-1}} F^*L\ast{F^m}^*M\xto{\psi_L\ast \psi_M}L\ast M
\eeq
in $\D_G^{F^{m+1}}(G)$ where $\psi_L$ and $\psi_M$ are the equivariance isomorphisms associated with $L\in \D^{F^m}_G(G)^F$ and $M\in \D^F_G(G)^{F^m}$ respectively. Recall from \cite[\S2.4.5]{De}, that the category $\D^{F^{m+1}}_G(G)$ has a normalized trace denoted by $\tr_{F^{m+1}}$ in the notation from {\it loc. cit.} The following is a generalization of \cite[Thm. 2.14]{De}:
\bthm\label{ipf}
Let $M\in \D^F_G(G)$ and $L\in \D_G^{F^{m}}(G)^F$. Then we have the functions ${N_m^{-1}}^*T_{L,\psi_L}$ and $\chi_M\in \FunG$. Then we have the following equality
\beq
\<{N_m^{-1}}^*T_{L,\psi_L},\chi_M\>=\tr_{F^{m+1}}(\zeta_{L,\eta_m(M)}).
\eeq
\ethm

\bcor\label{c:ipf}
Let $V\in \Irrep(G,F)$ with corresponding local system $V_{loc}\in \D^F_G(G)$. Let $W\in \Irrep(G,F^m)^F$ and let $(W_{loc},\psi_W)\in D^{F^m}_G(G)^F$ where $\psi_W$ is chosen as in \S\ref{s:irrepfrob}. Then we have the following relationship between the $m$-th Shintani basis and $\Irrep(G,F)$:
\beq
\<\Sh_m(W),\chi_V\>=\tr_{F^{m+1}}(\zeta_{(W_{loc},\psi_W),\eta_m(V_{loc})}).
\eeq
\ecor

\subsection{Proof of Theorem \ref{ipf}}\label{s:poipf}
The proof of Theorem \ref{ipf} is similar to the proof of \cite[Thm. 2.14]{De} as well as the proof of Theorem \ref{t:twisttrace} above. 

Let $\eta_m(M)=(M,\psi_M)$. We will often abuse notation and denote this object simply as $M\in \D^F_G(G)^{F^m}$. As before let $\phi_L,\phi_M$ denote the equivariance structures on $L,M$ corresponding to the twisted conjugations actions. Let us compute the stalk of the automorphism $\zeta_{L,M}$ at a point $t\in G$. By the definition of the crossed braidings, the stalk
$$\zeta_{L,M}(t):(L\ast M)_t\rar{} (L\ast M)_t$$ is given by the composition
$$(L\ast M)_t=\int\limits_{h_1F^m(h_2)=t} L_{h_1}\otimes M_{h_2} $$
$$\xto{{\int}\phi_L(h_2^{-1},h_1)\otimes \id_{M}}\int\limits_{h_1F^m(h_2)=t} L_{h_2^{-1}h_1F^m(h_2)}\otimes M_{h_2}=(M\ast F^*L)_t$$
$$\xto{\int\id_{L}\otimes \phi_M(F^{-1}(F^m(h_2)^{-1}h_1^{-1}h_2),h_2)}\int\limits_{h_1F^m(h_2)=t} F^*L_{F^{-1}(h_2^{-1}h_1F^m(h_2))}\otimes M_{F^{-1}(F^m(h_2)^{-1}h_1^{-1}h_2)h_1F^m(h_2)}$$
$$=(F^*L\ast{F^m}^*M)_t$$
$$\xto{\int\psi_{L}(F^{-1}(h_2^{-1}h_1F^m(h_2)))\otimes \psi_M(F^{-1}(h_2^{-1}F^{-m}(h_1^{-1}h_2))F^{-m}(h_1)h_2)}$$
$$\int\limits_{h_1F^m(h_2)=t} L_{F^{-1}(h_2^{-1}h_1F^m(h_2))}\otimes M_{F^{-1}(h_2^{-1}F^{-m}(h_1^{-1}h_2))F^{-m}(h_1)h_2}$$
%$$\xto{{\displaystyle\int}\left(\psi_L(F^{-1}(F^{-m}(h_2)^{-1}h_1h_2))\circ\phi_L(F^{-m}(h_2)^{-1},h_1)\right)\otimes \left(\psi_M(F^{-1}(h_2^{-1}h_1^{-1}F^{-m}(h_2))h_1h_2)\circ\phi_M(F^{-1}(F^m(h_2^{-1}h_1^{-1})h_2), h_2)\right)}$$
$$=(L\ast M)_t.$$

Consider the automorphism of the antidiagonal $\bar{\Delta}_t:=\{(h_1,h_2)\in G\times G|h_1F^m(h_2)=t\}$ that appears implicitly in the stalk of the automorphism $\zeta_{L,M}$ as above, namely
$$(h_1,h_2)\mapsto (F^{-1}(h_2^{-1}h_1F^m(h_2)), F^{-1}(h_2^{-1}F^{-m}(h_1^{-1}h_2))F^{-m}(h_1)h_2).$$
%$$=({F^{-1}(t^{-1}h_1t)}, {F^{-1}(t^{-1}h^{-1}_1t)}t).$$ 

The fixed point set of this automorphism is the finite set $\{(h_1,h_2)\in \bar{\Delta}_t|(h_1,h_2)\in R_{F^m,F}\}$. If we identify $\bar{\Delta}_t$ and $G$ using the first projection, then the inverse of this automorphism corresponds to the Frobenius automorphism of $G$ (as a scheme) defined by $h_1\mapsto tF^{m+1}(h_1)F^m(t)^{-1}$. Hence by the Grothendieck-Lefschetz trace formula and using the definition of the trace functions $T_{L,\psi_L},T_{M,\psi_M}$, we deduce that 
\beq
\tr(\zeta_{L,M}(t))=\sum\limits_{\substack{(h_1,h_2)\in R_{F^m,F}\\h_1F^m(h_2)=t}}T_{L,\psi_L}(h_1,h_2)\cdot T_{M,\psi_M}(h_2,h_1).
\eeq

Now $L\ast M$ is an object of $\D_G^{F^{m+1}}(G)\cong \bigoplus\limits_{\substack{<t> \in \\ H^1(F^{m+1},G)}}D^b\h{Rep}(G^{tF^{m+1}})$. Hence by definition of the trace $\tr_{F^{m+1}}$ on $\D^{F^{m+1}}_G(G)$ (see \cite[\S2.4.5]{De}) we have 
\beq
\tr_{F^{m+1}}(\zeta_{L,M})=\sum\limits_{\substack{<t> \in \\ H^1(F^{m+1},G)}}\frac{\tr(\zeta_{L,M}(t))}{|G^{tF^{m+1}}|}
\eeq
\beq
=\sum\limits_{\substack{<t> \in \\ H^1(F^{m+1},G)}}\sum\limits_{\substack{(h_1,h_2)\in R_{F^m,F}\\h_1F^m(h_2)=t}}\frac{T_{L,\psi_L}(h_1,h_2)\cdot T_{M,\psi_M}(h_2,h_1))}{|G^{tF^{m+1}}|}
\eeq
\beq
=\sum\limits_{\substack{<t> \in \\ H^1(F^{m+1},G)}}\sum\limits_{\substack{(h,t)=(h_1,h_1F^m(h_2))\in\\ R_{F^{m},F^{m+1}}}}\frac{T_{L,\psi_L}(h,F^{-m}(h^{-1}t))\cdot T_{M,\psi_M}(F^{-m}(h^{-1}t),h)}{|G^{tF^{m+1}}|}
\eeq
Note that in the last equality, we used the $G$-equivariant isomorphism $t_2:R_{F^m,F}\rar{\cong}R_{F^m,F^{m+1}}$ defined by (\ref{t2def}). Hence continuing further, we get $\tr_F(\zeta_{L,M})=$
\beq
=\sum\limits_{\substack{<t> \in \\ H^1(F^{m+1},G)}}\sum\limits_{\substack{(h,t)\in\\ R_{F^{m},F^{m+1}}}}\frac{T_{L,\psi_L}(h,F^{-m}(h^{-1}t))\cdot T_{M,\psi_M}(F^{-m}(h^{-1}t),h)}{|\Stab_G(h,t)|\cdot |{G^{tF^{m+1}}}\mbox{-orbit of }(h,t)|}
\eeq
\beq
=\sum\limits_{\substack{<(h,t)>\in G\backslash R_{F^{m},F^{m+1}}}}\frac{T_{L,\psi_L}(h,F^{-m}(h^{-1}t))\cdot T_{M,\psi_M}(F^{-m}(h^{-1}t),h)}{|\Stab_G(h,t)|}
\eeq
\beq
=\sum\limits_{\substack{<(h_1,h_2)>\in G\backslash R_{F^{m},F}}}\frac{T_{L,\psi_L}(h_1,h_2)\cdot T_{M,\psi_M}(h_2,h_1)}{|\Stab_G(h_1,h_2)|}
\eeq
\beq
=\sum\limits_{\substack{<(g,h)>\in G\backslash R_{\id_G,F}}}\frac{{N_m^{-1}}^*T_{L,\psi_L}(g,h)\cdot {N_m^{-1}}^*\tau^*T_{M,\psi_M}(g,h)}{|\Stab_G(g,h)|}
\eeq
\beq\label{eqchi}
=\sum\limits_{\substack{<(g,h)>\in G\backslash R_{\id_G,F}}}\frac{{N_m^{-1}}^*T_{L,\psi_L}(g,h)\cdot \chi_{M}(g^{-1},h)}{|\Stab_G(g,h)|}
\eeq
\beq
=<{N_m^{-1}}^*T_{L,\psi_L},\chi_M>
\eeq
as desired. Here in (\ref{eqchi}),  we have used Lemma \ref{twistandsf} and Example \ref{charfun}.

\section{Shintani descent for neutrally unipotent groups}\label{s:pfmain}
In the remainder of this paper, we restrict our attention to neutrally unipotent groups $G$ over $\k$ equipped with an $\Fq$-Frobenius $F$ . In this section, we will describe a refined version (see Theorem \ref{t:refthm}) of Theorem \ref{mainthm} and deduce that the theorem in fact follows from this refinement. We will prove Theorem \ref{t:refthm} in \S\ref{s:pfref}. We will use the theory of character sheaves on neutrally unipotent groups. Let us begin by recalling some aspects of this theory. 

\subsection{$\f{L}$-packets of characters of neutrally unipotent groups}\label{lpofchar}
In  \cite{B} a partitioning of the set $\Irrep(G,F)$ of irreducible characters  into $\f{L}$-packets  has been constructed. Let us recall this notion. Let $\widehat{G}$ denote the set of isomorphism classes of minimal idempotents in the braided monoidal category $\DG$. The $\f{L}$-packets of irreducible characters are parameterized by the set $\widehat{G}^F$ of $F$-stable minimal idempotents in  $\DG$. 
\bdefn\label{d:lp}
Let $e\in \widehat{G}^F$ and let $\psi_e:F^*e\rar\cong e$ be the unique isomorphism such that $T_e:=T_{e,\psi_e}\in \FunG$ is an idempotent. Let $W\in \Irrep(G,F)$. We say that $W$ lies in the $\f{L}$-packet associated with $e$, or equivalently that $W\in \Irrep_e(G,F)$ if and only if the following equivalent  (see \cite[Thm. 2.27]{De}) conditions hold
\bit
\item[(i)] The character $\chi_W$ lies in $T_e\FunG\subset \FunG$.
\item[(ii)] The idempotent $T_e$ acts by the identity on $W$. To be more precise, suppose that $W$ is an irreducible representation of a pure inner form $G^{gF}$. Then we want that the idempotent $T_e(\cdot,g)\in \Fun(G^{gF})^{G^{gF}}$ acts trivially on $W$.
\item[(iii)] The local system $W_{loc}\in e\D^F_G(G)\subset \D^F_G(G).$
\eit

\edefn

In \cite{BD}, character sheaves and their $\f{L}$-packets  are also defined. The $\f{L}$-packets of character sheaves are parameterized by $\widehat{G}$. Let us recall this:
\bdefn
Let $e\in \DG$ be a minimal idempotent. Let $\Meg\subset e\DG$ be the modular category associated with $e$ as defined in \cite{BD}. Then the set $CS_e(G)$ of character sheaves in the $\f{L}$-packet associated with $e$ is defined to be the set of simple objects (up to isomorphism) of the modular category $\Meg\subset \DG$. The set $CS(G)$ of all character sheaves on $G$ is defined to be the union (which turns out to be disjoint) of all the sets $CS_e(G)$ as $e$ ranges over $\widehat{G}$.
\edefn

We recall the following results from \cite{B}:
\bthm\label{t:bmain}
(i) We have a partition
\beq
\Irrep(G,F)=\coprod\limits_{e\in \widehat{G}^F}{\Irrep_e(G,F)}
\eeq and an isomorphism of algebras
\beq
\FunG=\bigoplus\limits_{e\in \widehat{G}^F}T_e\FunG.
\eeq
(ii) The set $\Irrep_e(G,F)$ is (after taking characters of the irreducibles) an orthonormal basis of $T_e\FunG$.\\
(iii) For each $C\in CS_e(G)^F$, we choose $\psi_C:F^*C\rar{\cong} C$ so that \footnote{It is possible to impose some additional conditions such that $\psi_C$ will be uniquely determined up to scaling by roots of unity.} $|T_{C,\psi_C}|=1$. Then the set $\{T_{C,\psi_C}\}_{C\in CS_e(G)^F}$ is also an orthonormal basis of $T_e\FunG$.
\ethm

\subsection{Shintani descent and $\f{L}$-packets}
In this section we show that Shintani descent respects the $\f{L}$-packet decomposition. In particular, we will show that to prove Theorem \ref{mainthm}, it is enough to restrict our attention to each $\f{L}$-packet.

For any positive integer $m$ we have 
\beq
\Irrep(G,F^m)=\coprod\limits_{e\in \widehat{G}^{F^m}}\Irrep_e(G,F^m).
\eeq
Moreover by Definition \ref{d:lp}(iii) and the fact that $F^*$ preserves convolution, the permutation of $\Irrep(G,F^m)$ induced by $F$ (see \S\ref{s:agof}, \S\ref{s:irrepfrob}) preserves this $\f{L}$-packet decomposition:
\beq
F^*:\Irrep_e(G,F^m)\rar{\cong}\Irrep_{F^*e}(G,F^m).
\eeq
Hence
\beq\label{e:lpd}
\Irrep(G,F^m)^F=\coprod\limits_{e\in \widehat{G}^{F}}\Irrep_e(G,F^m)^F.
\eeq

Recall that Shintani descent is a map $\Sh_m:\Irrep(G,F^m)^F\hookrightarrow \FunG$.

\bprop
Let $e\in \widehat{G}^F$. Then we have $\Sh_m(\Irrep_e(G,F^m)^F)\subset T_e\FunG$ and that the image is an orthonormal basis of $T_e\FunG$. In other words Shintani descent preserves $\f{L}$-packets.
\eprop
\bpf
Let $W\in \Irrep_e(G,F^m)^F$ and let $\psi_W$ be an intertwiner as in \S\ref{s:irrepfrob}. Since $W$ lies in the $\f{L}$-packet associated with $e$, $W_{loc}\in e\D^{F^m}_G(G)$, i.e.  $e\ast W_{loc}\cong W_{loc}$. Hence by Proposition \ref{s:convandsf}, $T_{W_{loc},\psi_W}\in T_e\Fun_G(R_{F^m,F})\subset \Fun_G(R_{F^m,F}).$ Now by Lemma \ref{l:modiso}, ${N_m^{-1}}^*:\Fun_G(R_{F^m,F})\to \FunG$ is a $\FunG$-module isomorphism. Hence $\Sh_m(W)={N_m^{-1}}^*T_{W_{loc},\psi_W}\in T_e\FunG$ as desired. Moreover, since $\Sh_m(\Irrep(G,F^m)^F)$ is an orthonormal basis of $\FunG$, by Theorem \ref{t:bmain}(i) and (\ref{e:lpd}) $\Sh_m(\Irrep_e(G,F^m)^F)\subset T_e\FunG$ must in fact be an orthonormal basis.
\epf

\brk\label{r:shme}
We say that $\Sh_m(\Irrep_e(G,F^m)^F)$ is the $m$-th Shintani basis of the subspace $T_e\FunG$ associated with the minimal idempotent $e$.
\erk

\subsection{A refined version of Theorem \ref{mainthm}}
We can now state the refined version of our main result, Theorem \ref{mainthm} and deduce the main theorem from this refinement.
\bthm\label{t:refthm}
Let $e\in\widehat{G}^F$ be an $F$-stable minimal idempotent. Then we have the following:\\
(i) There exists a positive integer $m_0$ such that if $m$ is divisible by $m_0$, then the $m$-th Shintani basis of $T_e\FunG$ (see Remark \ref{r:shme}) is (up to scaling by roots of unity) independent of $m$. We define almost characters associated with $e$ to be the elements of this common Shintani basis.\\
(ii) The almost characters associated with $e$ as defined above are eigenvectors for the twisting operator $\Theta$ (see (\ref{twistintro})).\\
(iii) The mapping which takes an $F$-stable character sheaf $C\in CS_e(G)^F$ to its associated trace of Frobenius function $T_{C,\psi_C}$ in $T_e\FunG$ defines a bijection from the set $CS_e(G)^F$ to the set of almost characters associated with $e$.
\ethm

 We will prove the theorem in \S\ref{s:pfref} below.

To prove Theorem \ref{mainthm}(i) from this, it suffices to take a common multiple of all the $m_0$'s corresponding to all the $e$'s in the finite set $\widehat{G}^F$. Parts (ii), (iii) of Theorem \ref{mainthm} also follow readily from the above refinement.

\section{Categorical preliminaries}\label{s:catprel}

We have the $F$-stable minimal idempotent $e$ in $\DG$. Note that we have the integer $n_e$ from \cite[Thm. 1.15]{BD} such that the category $\Meg[-n_e]\subset \DG$ consists of perverse sheaves. 

Let $\M_{GF,e}\subset e\D_G(GF)\cong e\D^F_G(G)$ be the full subcategory formed by those objects whose underlying $\Qlcl$-complex is a perverse sheaf shifted $n_e$. By \cite[Thm. 2.27]{De}, $\M_{GF,e}$ is an invertible $\Meg$-module category. By the results of \cite{BD},\cite{De} we also have $e\DG\cong D^b\Meg$ and $e\D_G^{F^m}(G)\cong D^b\M_{GF^m,e}$.

To complete the proof of Theorem \ref{t:refthm}, we will use some categorical notions and results from \cite{ENO}. In this section we recall some of these facts and apply them in our setting of neutrally unipotent groups. 

\subsection{Twists in module categories}\label{s:twistinmod}
By \cite[Thm. 5.2]{ENO}, we have an equivalence of groupoids 
\beq
\uPic(\Meg)\cong \uEqBr(\Meg).
\eeq
Under this equivalence the invertible $\Meg$-module category $\M_{GF,e}$ corresponds to the modular autoequivalence $F={F^{-1}}^*:\Meg \rar{\cong} \Meg$ since for $C\in \Meg\subset e\DG$ and $M\in \M_{GF,e}\in e\D^F_G(G)$ we have the functorial crossed braidings $\beta_{M,C}:M\ast C\to F(C)\ast M$.

Let us consider the semidirect product $G\<F\>=G\rtimes \Z$. The categories 
\beq
e\D_G(G\<F\>):=\bigoplus\limits_{m\in \Z}e\D_G(GF^m) \mbox{ and}
\eeq 
\beq
\M_{G\<F\>,e}:=\bigoplus\limits_{m\in \Z}\M_{GF^m,e}
\eeq 
are (infinite) spherical braided $\Z$-crossed category with trivial components $e\DG$ and $\Meg$ respectively with a (rigid) duality functor which we denote by ${(\cdot)}^\vee$ (cf. \cite[\S2.3]{De1}) and a natural  isomorphism $(\cdot)^{\vee\vee}\cong \id$ of monoidal functors. In particular, for each $M\in e\D^F_G(G)$ we have a twist $\theta^{F,e}_M$ defined as the composition 
\beq\label{twistinmod}
\theta^{F,e}_M: M \xto{\id_M\ast \h{coev}_M} M\ast M\ast M^\vee \xto{\beta_{M,M}} F(M)\ast M\ast M^\vee\xto{\id_{F(M)}\ast{\h{ev}_{M^\vee}}}F(M).
\eeq

From \cite{BD}, \cite{De1}, it is known that the twist $\theta:\id_{\DG}\to \id_{\DG}$ in $\DG$ defines the spherical structure of the modular category $\Meg\subset e\DG$.  Similarly we can define twists $\theta^{F^m,e}$ in each $e\DG$-module category $e\D^{F^m}_G(G)$. Then using the properties of spherical braided crossed categories (for example using string diagrams to aid visualization), we can prove an analogues of Lemma \ref{modcattwist} and Theorem \ref{t:twisttrace}:
\blem\label{l:modulartwisttrace}
(i) Let $C\in e\DG$ and let $M\in e\D^F_G(G)$. Then $\theta^{F,e}_{C\ast M}$ is equal to the composition
\beq
C\ast M\xto{\beta_{C,M}} M\ast C\xto{\beta_{M,C}} F(C)\ast M \xto{\theta_{F(C)}\ast \theta^{F,e}_M}F(C)\ast F(M)\rar{} F(C\ast M).
\eeq
(ii) In fact, all the twists $\theta^{F^a,e}$ in $\D^{F^a}_G(G)$ for $a\in \Z$ are compatible with each other, namely if $L\in \D^{F^a}_G(G)$, $M\in \D^{F^b}_G(G)$ then $\theta^{F^{a+b},e}_{L\ast M}$ equals the composition
\beq
L\ast M\xto{\beta_{L,M}} F^a(M)\ast L\xto{\beta_{F^a(M),L}} F^b(L)\ast F^a(M) \xto{\theta^{F^a,e}_{F^b(L)}\ast \theta^{F^b,e}_{F^a(M)}}F^{a+b}(L)\ast F^{a+b}(M)\rar{} F^{a+b}(L\ast M).
\eeq
(iii) For $M\in e\D^F_G(G)$ we define the automorphism $\nu_{\theta^{F,e}_M}$ in $e\D^{F^2}_G(G)$ to be the composition
\beq
\nu_{\theta^{F,e}_M}:M\ast M\xto{\beta_{M,F^*M}^{-1}} M\ast F^*M\xto{\id_M\ast F^*(\theta^{F,e}_M)}M\ast M.
\eeq
Then we have
\beq
\tr_{F,e}(\id_M)=\tr_{F^2,e}(\nu_{\theta^{F,e}_M}),
\eeq
where $\tr_{F,e}$ denotes the $e\DG$-module trace (cf. \cite{Sch}) in $e\D^F_G(G)$ coming from the natural spherical structure.
\elem

Combining this with Lemma \ref{modcattwist}, we conclude that both the twists $\theta^F$ and $\theta^{F,e}$ in the invertible $\Meg$-module category $\M_{GF,e}$ are compatible with the twist $\theta$ in $\Meg$. Hence we must have
\bprop\label{twotwists}
For each $m\in \Z$, there is a constant $c_{F^m,e}\in \Qlcl^\times$ such that for each $M\in e\D^{F^m}_G(G)$, we have $\theta^{F^m}_M=c_{F^m,e}\cdot \theta^{F^m,e}_M$ and that $c_{F^m,e}=(c_{F,m})^m$ for each $m\in \Z$. In fact we must have 
\beq\label{e:wantcfe}
c_{F,e}=\sigma_{F,e}\cdot \frac{q^{\dim G}}{q^{d_e}},
\eeq
where $d_e:=\frac{\dim G-n_e}{2}$ is called the functional dimension of $e$ and where $\sigma_{F,e}=\pm 1$.
\eprop
\bpf
The equality (\ref{e:wantcfe}) is proved by (\ref{e:cfe}) below. The remaining part of the proposition is clear from the previous results and the uniqueness (cf. \cite{Sch}) of module traces up to scaling.
\epf

\subsection{Positivity of spherical structure}
By \cite[Thm. 2.17]{De2}, the Frobenius-Perron dimension of $\Meg$ is an integer. In this section, let us choose an identification $\Qlcl\cong \f{C}$. Then by \cite[Cor. 2.24]{DGNO} the braided monoidal category $\M_{G,e}$ has a unique positive (with respect to our chosen identification) spherical structure. Let us denote the corresponding twist by $\theta^+:\id_{\Meg}\to \id_{\Meg}$. Similarly the crossed braided category $\M_{G\<F\>,e}$ also has a unique positive spherical structure. Using this we obtain new twists $\theta^{F^m,e,+}$ in the module categories $\M_{GF^m,e}$ for each $m\in \Z$. Let $\M_{G\<F\>,e}^+$ denote the braided $\Z$-crossed category equipped with this positive spherical structure. For each $m\in \Z$, let $\tr_{F^m,e}^+$ denote the $\Meg^+$-module trace in $\M^+_{GF^m,e}$ corresponding to the positive spherical structure. Using \cite[\S2.4.3]{DGNO} we deduce
\blem\label{l:sgnm}
For $m\in \Z$ and a simple object $M\in \M_{GF^m,e}$, we must have $$\theta_M^{F^m,e,+}=\sigma_M\theta^{F^m,e}_M:M\rar{\cong} F^m(M)$$
where $\sigma_M=\pm 1$.
\elem

\brk
Since $\theta^+, \theta^{F^m,e,+}, \tr^+_{F^m,e}$ are all defined using a spherical structure on a braided $\Z$-crossed category, they must satisfy analogues of Lemma \ref{l:modulartwisttrace}. In particular, for each $M\in \M_{GF,e}^+$ we must have 
\beq\label{e:postr}
\tr_{F,e}^+(\id_M)=\tr_{F^2,e}^+(\nu_{\theta^{F,e,+}_M}).
\eeq
\erk

Let $M\in \M_{GF,e}$ be simple object. Now by Theorem \ref{t:twisttrace}, we have
\beq
\tr_{F}(\id_M)=\tr_{F^2}(\nu_{\theta^{F}_M})=\tr_{F^2}(\nu_{\sigma_M\cdot c_{F,e}\cdot\theta^{F,e,+}_M})=\sigma_M\cdot c_{F,e}\cdot \tr_{F^2}(\nu_{\theta^{F,e,+}_M}).
\eeq
In \cite[Thm. 2.30]{De}, we have described the explicit relation between the trace $\tr_F$ and $\tr_{F,e}^+$. Using this we obtain
\beq
\frac{(-1)^{2d_e}}{\sqrt{\dim M_{G,e}}}\cdot\frac{q^{d_e}}{q^{\dim G}}\cdot\tr_{F,e}^+(\id_M)=\frac{(-1)^{2d_e}}{\sqrt{\dim M_{G,e}}}\cdot\frac{q^{2d_e}}{q^{2\dim G}}\cdot\sigma_M\cdot c_{F,e}\cdot \tr_{F^2,e}^+(\nu_{\theta^{F,e,+}_M}),
\eeq
where  $d_e$ is the functional dimension of $e$ (see Proposition \ref{twotwists}).
Comparing with (\ref{e:postr}), we obtain 
\beq\label{e:cfe}
c_{F,e}=\sigma_M\cdot\frac{q^{\dim G}}{q^{d_e}}.
\eeq
In particular this means that $\sigma_M$ does not depend on $M$ and we set $\sigma_{F,e}=\sigma_M$. This implies the natural trace $\Meg$-module trace $\tr_{F,e}$ in $\M_{GF,e}$ is either purely positive or purely negative. Hence the natural spherical structure on $\Meg$ must in fact be positive. Moreover, this is true with respect to any identification $\Qlcl\cong \f{C}$. This means that the weakly integral modular category $\Meg$ must actually be integral.
Hence we obtain the following result which confirms a conjecture from \cite{De}:
\bthm\label{t:posint}
Let $G$ be any neutrally unipotent group defined over $\k$ and let $e\in\DG$ be any minimal idempotent. Then the modular category $\M_{G,e}$ is positive integral.
\ethm
\bpf
By the results of \cite{BD}, the minimal idempotent $e$ can be obtained from some admissible pair\footnote{We refer to {\it op. cit.} for a precise definition.} $(H,\L)$. Now the subgroup $H$ must be stable under some Frobenius map $F:G\to G$. In fact we may choose a Frobenius $F$ so that the multiplicative local system $\L\in H^*$ (the Serre dual of $H$, cf. {\it op. cit.}) is also $F$-stable. Then the minimal idempotent $e$ is also $F$-stable for such an $F$. Now we are in the situation where we can apply our previous results to deduce the theorem. 
\epf

Using the observations of this section and \cite[Thm. 2.30]{De}, we see that:
\bcor\label{c:comptrace}
There is the following relationship between the $e\DG$-module traces $\tr_F$ and $\tr_{F,e}$ in the category $e\D^F_G(G)$
\beq
\tr_{F,e}=\sigma_{F,e}\cdot(-1)^{2d_e}{\sqrt{\dim \Meg}}\cdot\frac{q^{\dim G}}{q^{q^{d_e}}}\cdot \tr_F=(-1)^{2d_e}{\sqrt{\dim \Meg}}\cdot c_{F,e}\cdot \tr_F.
\eeq
\ecor

\subsection{Braided monoidal actions and braided crossed categories}

Now by \cite[Thm. 4.5]{ENO}, we must have $F^{n'}\cong \id_{\Meg}$ as braided monoidal functors for some positive integer $n'$. This gives us a group homomorphism ${\Z}/{n'\Z}\to \Pic(\Meg)\cong \EqBr(\Meg)$. However, this may not necessarily give rise to a braided monoidal action of $\Z/n'\Z$ on $\Meg$. By \cite[\S8]{ENO}, the obstruction is described by a 3-cocycle in $H^3(\Z/n'\Z, Z)$ where $Z$ is the group of automorphism of $\id_{\Meg}$ considered as a braided tensor functor or equivalently $Z=\Aut_{\Meg\h{-Mod}}(\Meg)$. We can trivialize this 3-cocycle by choosing a suitable multiple $n$ of $n'$. This means that we have a natural equivalence $\id_{\Meg}\rar{\xi}F^n$ (of braided monoidal functors) such that this defines a braided monoidal action of $\Z/n\Z$ on $\Meg$ or equivalently a map $\Z/n\Z\rar{}\uEqBr(\Meg)\cong \uPic(\Meg)$. Now note that in fact the cohomology group $H^4(\Z/n\Z,\Qlcl^\times)$ is trivial. Hence by \cite[\S8]{ENO} we can in fact lift this to obtain a map $\Z/n\Z\rar{}\uuPic(\Meg)$ i.e. we can put a structure of a braided $\Z/n\Z$-crossed category on the category 
\beq
\D:=\bigoplus\limits_{i=0}^{n-1}\M_{GF^i,e}.
\eeq
This means that we have an induced action of $\Z/n\Z$ on each of the categories $\M_{GF^a,e}$ compatible with the braided monoidal action of $\Z/n\Z$ on $\Meg$. We have well defined natural isomorphisms $F^a\rar{\xi^m_a}F^{a+m}$ for each $a\in \Z$ and $m\in n\Z$.  Also we have the corresponding equivalences of $\Meg$-module categories
\beq
\t\xi^m_a:\M_{GF^a,e}\rar{\cong}\M_{GF^{a+m},e}
\eeq
for $a\in \Z, m\in n\Z$ and these satisfy certain compatibility conditions. The functor $\t\xi_a^m$ commutes with the actions of $\Z/n\Z$. 

Note that we have the equivalence of $\Meg$-module categories
\beq
\t\xi=\t\xi_0^n: \M_{G,e}\rar{\cong}\M_{GF^{n},e}.
\eeq
In particular we have the object $\t\xi e\in \M_{GF^n,e}$ and the functor $\t\xi$ above is canonically equivalent to the functors $C\mapsto \t\xi e\ast C$ and $C\mapsto C\ast \t\xi e$. And similarly for each $a\in \Z, m\in n\Z$ and $M\in \M_{GF^a,e}$, we have natural isomorphisms $(\t\xi e)^{\ast\frac{m}{n}}\ast M\cong\t\xi_a^m(M)\cong M\ast (\t\xi e)^{\ast\frac{m}{n}}\rar{\cong} F^a(\t\xi e)^{\ast\frac{m}{n}}\ast M$ giving us a canonical isomorphism $\t\xi e \cong F^a(\t\xi e)$ for each $a\in \Z$. 

\brk
The fact that $\bigoplus\limits_{i\in \Z}\M_{GF^i,e}$ is a spherical braided $\Z$-crossed category implies that each $\M_{GF^i,e}$ is equipped with a normalized $\Meg$-module trace (in the sense of \cite[\S2.3.1 (2)]{De}). Hence under the $\Meg$-module equivalence $\t\xi^n_0$, the $\Meg$-module traces on $\Meg$ and $\M_{GF^n,e}$ must either agree or be negatives of each other. In anycase the traces on $\Meg$ and $\M_{GF^{2n},e}$ must agree. Hence replacing $n$ by $2n$ if required, we may assume that the traces agree and hence we obtain a spherical structure on the category $\D$.
\erk

Note that all the structure above can be extended to the corresponding bounded derived categories. Now recall from (\ref{defzeta}) that for $C\in e\DG^F$ and $M\in e\D^F_G(G)$ we have the composition
\beq
\zeta_{C,M}: C\ast M\xto{\beta_{M,F^*C}^{-1}} M\ast F^*C \xto{\beta_{F^*C,M}^{-1}} F^*C\ast M\xto{\psi_C\ast \id_M}C\ast M
\eeq
in the category $e\D^F_G(G).$ Then we have
\blem\label{l:xizeta}
Let $m\in n\Z$. Then the automorphism $\t\xi_1^m(\zeta_{C,M})$ in $e\D^{F^{m+1}}_G(G)$ is equal to the composition
{\beq
\t\xi_1^m(\zeta_{C,M}): \t\xi_0^mC\ast M\xto{\beta^{-1}} M\ast F^*\t\xi_0^mC\xto{\beta^{-1}} F^*\t\xi_0^mC\ast {F^m}^*M\xto{\t\xi^m_0(\psi_C)\ast ({\xi_0^m})_{{F^m}^*M}}\t\xi^m_0C\ast M.
\eeq}
\elem

\subsection{Certain twists are roots of unity}
Recall from (\ref{twistinmod}) that we have the twist $\theta^{F,e}$ in the $\Meg$-module category $\M_{GF,e}$. For each $M\in \M_{GF,e}$ we have the following automorphism of $M$
\beq\label{thetan}
\omega_M:M\xto{\theta^{F,e}_M}F(M)\xto{F(\theta^{F,e}_M)}\cdots\xto{F^{n-1}(\theta^{F,e}_M)}F^n(M)\xto{(\xi_0^n)_M^{-1}}M. 
\eeq

\blem\label{rootofunity}
If $M\in \M_{GF,e}$ is a simple object then $\omega_M$ is a root of unity.
\elem
\bpf
We have the spherical braided $\Z/n\Z$-crossed category $\D$. Taking the equivariantization we obtain the modular category $\D^{\Z/n\Z}$. We choose an isomorphism $M\xto{\psi_M}F(M)$ such that $(M,\psi_M)$ defines a $\Z/n\Z$-equivariant object in the category $\M_{GF,e}^{\Z/n\Z}\subset \D^{\Z/n\Z}$. Or in other words we must have
\beq\label{psin}
\id_M:M\xto{\psi_M}F(M)\xto{F(\psi_M)}\cdots\xto{F^{n-1}(\psi_M)}F^n(M)\xto{(\xi_0^n)_M^{-1}}M. 
\eeq 
Now since $\D^{\Z/n\Z}$ is a modular category, it is equipped with a twist $\theta'$. The twist $\theta'_{M,\psi_M}$ is equal to the composition $M\xto{\theta^{F,e}_M}F(M)\xto{\psi_M^{-1}}M$. Now we know that twists of simple objects in a modular category are roots of unity. Hence $\theta'_{M,\psi_M}$ is a root of unity. On the other hand, by comparing (\ref{thetan}) and (\ref{psin}) we see that $\omega_M=({\theta'}_{M,\psi_M})^n$. Hence $\omega_M$ is also a root of unity.
\epf

\bcor\label{d:m0}
There exists a positive integer $m_0$ divisible by $n$ such that for each $m\in m_0\Z\subset n\Z$ and each simple object $M\in \M_{GF,e}$ the composition
\beq
M\xto{\theta^{F,e}_M}F(M)\xto{F(\theta^{F,e}_M)}\cdots\xto{F^{m-1}(\theta^{F,e}_M)}F^{m}(M)\xto{(\xi_0^{m})_M^{-1}}M
\eeq
is equal to the identity or equivalently (see Lemma \ref{twotwists}), such that the composition
\beq\label{m0prop}
M\xto{\theta^{F}_M}F(M)\xto{F(\theta^{F}_M)}\cdots\xto{F^{m-1}(\theta^{F}_M)}F^{m}(M)\xto{(\xi_0^{m})_M^{-1}}M
\eeq
equals scaling by $(c_{F,e})^m$.\\
%(ii) Let $C$ be an $F$-stable simple and let $\psi_C:C\rar{\cong}F(C)$ be such that the composition
%\beq
%C\xto{\psi_C}F(C)\xto{F(\psi_C)}\cdots\xto{F^{n-1}(\psi_C)}F^n(C)\xto{(\xi_0^n)_C^{-1}}C
%\eeq
%equals $\id_C$.
\ecor

\bpf
Since $\M_{GF,e}$ has finitely many simple objects, there exists a positive integer $b$ such that $\omega_M^b=1$ for all simples $M$. Then it is easy to check that $m_0=nb$ satisfies the desired property. 
\epf

\section{Completion of the proof}\label{s:pfref}
We can now complete the proof of Theorem \ref{t:refthm}. Let $M\in \D^F_G(G)$ and let $m$ be any integer. Then we have the object $\eta_m(M)=(M,\psi_{M,m})\in \D^F_G(G)^{F^m}$ defined by Lemma \ref{d:eta}. Now suppose that $m\in m_0\Z$, where $m_0$ is as in Corollary \ref{d:m0}. Now applying ${F^m}^*$ to (\ref{m0prop}) and comparing with (\ref{e:eta}) we see that 
\beq\label{e:psixi}
\psi_{M,m}=(c_{F,e})^m\cdot {F^m}^*(\xi_0^m)_M=(c_{F,e})^m\cdot(\xi_0^m)_{{F^m}^*M}:{F^m}^*M\to M
\eeq for each $m\in m_0\Z$.

We will now prove that the integer $m_0$ above satisfies the desired properties of Theorem \ref{t:refthm}. Let $m$ be any positive multiple of $m_0$. Let $C\in \Meg$ be an $F$-stable simple object. We choose an $F$-equivariance isomorphism $\psi_{C,m}:F^*C\rar{\cong}C$ such that the composition
\beq\label{psicondi}
C\xto{F(\psi_{C,m})}F(C)\xto{F^2(\psi_{C,m})}\cdots\xto{F^{m}(\psi_{C,m})}F^m(C)\xto{(\xi_0^m)_C^{-1}}C
\eeq 
equals the twist $\theta_C:C\to C$. Now consider the $F$-stable object $L_m:=\t\xi_0^m(C)$ in $\M_{GF^m,e}\subset e\D^{F^m}_G(G)$. It is equipped with the associated $F$-equivariance isomorphism 
\beq
\psi_{L_m}:=c_{F,e}\cdot{\t\xi_0^m(\psi_{C,m})}:F^*L_m=F^*\t\xi_0^m(C)=\t\xi^m_0(F^*C)\rar{\cong}\t\xi_0^m(C)=L_m.
\eeq 
The reason for introducing the scaling factor $c_{F,e}$ in the definition of $\psi_{L_m}$ is that now the composition
\beq\label{psilcondi}
L_m\xto{F(\psi_{L_m})}F(L_m)\xto{F^2(\psi_{L_m})}\cdots\xto{F^{m}(\psi_{L_m})}F^m(L_m)
\eeq 
equals $(c_{F,e})^m\cdot \theta^{F^m,e}_{L_m}=c_{F^m,e}\cdot\theta^{F^m,e}_{L_m} =\theta^{F^m}_{L_m}$. Now $L_m[-\dim G-n_e]$ is an $F$-stable local system, say equal to ${W_m}_{loc}$ where $W_m\in \Irrep(G,F^{m})^F$. Then by (\ref{psilcondi}) the $F$-equivariance structure 
\beq
\psi_{W_m}:=\psi_{L_m}[-\dim G-n_e]:F^*{W_m}_{loc}\rar{\cong}{W_m}_{loc}
\eeq
satisfies the desired property as described in \S\ref{s:irrepfrob} and Remark \ref{r:altpsi}. In particular we must have
\beq\label{norm1}
||T_{W_m,\psi_{W_m}}||=||T_{L_m,\psi_{L_m}}||=||c_{F,e}\cdot T_{C,\psi_{C,m}}||=||T_{C,c_{F,e}\cdot \psi_{C,m}}||=1.
\eeq

Now combining Lemma \ref{l:xizeta} and (\ref{e:psixi}), we see that the composition (see (\ref{defzeta}))
\beq
\zeta_{(\t\xi_0^mC,\t\xi_0^m(\psi_{C,m})),\eta_mM}: \t\xi_0^mC\ast M\xto{\beta^{-1}} M\ast F^*\t\xi_0^mC\xto{\beta^{-1}} F^*\t\xi_0^mC\ast {F^m}^*M\xto{\t\xi^m_0(\psi_C)\ast \psi_{M,m}}\t\xi^m_0C\ast M
\eeq 
is equal to $(c_{F,e})^m\cdot \t\xi_1^m(\zeta_{(C,\psi_{C,m}),M})$. Hence we see that
\beq
\zeta_{(L_m,\psi_{L_m}),\eta_mM}=(c_{F,e})^{m+1}\cdot\t\xi^m_1(\zeta_{(C,\psi_{C,m}),M})
\eeq
and that
\beq
\zeta_{({W_m}_{loc},\psi_{W_m}),\eta_mM}=(-1)^{2d_e}\cdot(c_{F,e})^{m+1}\cdot\t\xi^m_1(\zeta_{(C,\psi_{C,m}),M})
\eeq

Now under the identification
\beq
\t\xi_1^m:e\D^F_G(G)\rar{\cong}e\D^{F^{m+1}}_G(G),
\eeq
the traces $\tr_{F,e}$ and $\tr_{F^{m+1},e}$ defined using the spherical structure on $e\D_G(G\<F\>)$ agree. Hence 
\beq
\tr_{F^{m+1},e}(\zeta_{({W_m}_{loc},\psi_{W_m}),\eta_mM})=(-1)^{2d_e}\cdot(c_{F,e})^{m+1}\cdot\tr_{F,e}(\zeta_{(C,\psi_{C,m}),M}).
\eeq
Hence using Corollary \ref{c:comptrace} we obtain
\beq
\tr_{F^{m+1}}(\zeta_{({W_m}_{loc},\psi_{W_m}),\eta_mM})=\frac{1}{\sqrt{\dim \Meg}}\cdot{\tr_{F,e}(\zeta_{(C,\psi_{C,m}),M})}
\eeq
\beq
=(-1)^{2d_e}c_{F,e}\cdot{\tr_{F}(\zeta_{(C,\psi_{C,m}),M})}=(-1)^{2d_e}{\tr_{F}(\zeta_{(C,c_{F,e}\cdot\psi_{C,m}),M})}.
\eeq
Now let $V\in \Irrep_e(G,F)$ and set $M=V_{loc}\in e\D^F_G(G)$. Then from Theorem \ref{ipf} and Corollary \ref{c:ipf} we deduce that
\beq
\<\Sh_m(W_m),\chi_V\>=(-1)^{2d_e}\<T_{C,c_{F,e}\cdot\psi_{C,m}},\chi_V\>.
\eeq
Now the twist $\theta_C$ is a root of unity. Hence the $F$-equivariance isomorphisms $\psi_{C,m}:F^*C\rar\cong C$ (chosen according to (\ref{psicondi})) for different $m\in m_0\Z$ differ from each other only up to scaling by roots of unity. Hence  $\Sh_m(W_m)$ for different values of $m$ (in $m_0\Z_{>0}$) as well as the trace of Frobenius functions $T_{C,c_{F,e}\cdot\psi_{C,m}}=c_{F,e}\cdot T_{C,\psi_{C,m}}$ all differ from each other only up to scalings by roots of unity. This completes the proof of Theorem \ref{t:refthm} (i) and (iii). Theorem \ref{t:refthm} (ii) also follows from Proposition \ref{p:sftwist}.

\end{document}